\documentclass[11pt,reqno]{amsart}
\usepackage{graphicx}
\usepackage{color}
\usepackage{mathrsfs,amsmath,amssymb,amsthm,amsfonts}

\usepackage{graphics,url,color,epsfig}
\usepackage{tikz}
\usepackage[colorlinks]{hyperref}
\AtBeginDocument{%
  \hypersetup{%
    linkcolor=blue,%
    citecolor=green,%
  }%
}
\usepackage{cleveref}
\usepackage{pgfplots}
\usepackage{multicol}

\def\R{\mathbb{R}}

\makeatletter

\newcommand{\Rmnum}[1]{\expandafter\@slowromancap\romannumeral #1@}
\makeatother

\newtheorem{thm}{Theorem}[section]

\newtheorem{lemma}[thm]{Lemma}
\newtheorem{remark}[thm]{Remark}
\newtheorem{theorem}[thm]{Theorem}

\usepackage[numbers,sort&compress]{natbib}
\allowdisplaybreaks

\begin{document}

\author{Liying Shan}
\address{Liying Shan \newline\indent School of Mathematics and Statistics \newline\indent Central China Normal University \newline\indent Wuhan 430079, P. R. China}
\email{shanliying@mails.ccnu.edu.cn}

\author{Wei Shuai}
\address{Wei Shuai \newline\indent
School of Mathematics and Statistics 
\newline\indent  Key Laboratory of Nonlinear Analysis and Applications
(Ministry of Education)
\newline\indent Central China Normal University \newline\indent Wuhan 430079, P. R. China}
\email{wshuai@ccnu.edu.cn}

\author{Leyun Wu}
\address{Leyun Wu
\newline\indent School of Mathematics \newline\indent South China University of Technology \newline\indent
Guangzhou, 510640, P. R. China}
\email{leyunwu@scut.edu.cn}

\title[ \tiny {Multiple solutions to the nonlinear Schr\"{o}dinger equation  with  a partial confinement}]{Multiple solutions to the nonlinear Schr\"{o}dinger equation with a partial confinement}

\begin{abstract}
We consider multiple solutions to the nonlinear Schr\"{o}dinger equation (NLS)  with  a partial confinement, which  is physically relevant to  dynamics of the Bose-Einstein condensate.
Our study not only verifies the existence of positive ground state solutions and the nonexistence of least energy sign-changing solutions but also sheds light on the symmetry associated with these solutions.

A novel finding is the existence of saddle type nodal solutions  with their nodal domains intersecting at the origin.
Furthermore, we have developed some innovative techniques such as the method of moving planes and the Hopf lemma for nonlinear Schr\"{o}dinger equations with partial confinement.
\end{abstract}



\maketitle

\numberwithin{equation}{section}

\section{Introduction}

In this paper, we focus on the standing waves of the type $\psi(z,t)=e^{-i\omega t}u(z)$ to the following
time dependent nonlinear Schr\"odinger equation
\begin{equation*}
i\partial_t\psi+ \Delta \psi-(x_1^2+x_2^2+\cdots+x_m^2+\omega)\psi+|\psi|^{p-2}\psi=0,\ \ (z,t)\in \R^N\times \R_+,
\end{equation*}
which is a basic quantum mechanics and nonlinear optics model,
where $i$ is the imaginary unit, $\omega$ is a real number, $N\geq 2$ and  $1\leq m \leq N$.

 The Cubic NLS for $p=4$, often referred to as the Gross-Pitaevskii equation (GPE), is of great importance in physics. This model specifically describes the mean-field dynamics of an extremely cold boson gas, which is a Bose-Einstein condensate confined in an anisotropic trap (see \cite{ACD,BBJV,F,GS}). We refer to Chen \cite{Chen-ARMA-2013} for a detailed derivation of the model \eqref{eqs1.1} from the many-body bosonic system.
 It is a well known fact that in this range the ground states of the translation invariant NLS, where the term
$(x_1^2+x_2^2+\cdots+x_m^2+\omega)\psi$ is removed, are unstable by blow up, as shown in \cite{BC-CRASP-1981, C1981}. However, the situation changes dramatically when we introduce a confinement term.

This leads to  study the existence of solutions to the following nonlinear Schr\"odinger equation (NLS)
 \begin{equation}\label{eqs1.1}
\begin{cases}
-\Delta u+(x_1^2+x_2^2+\cdots+x_m^2)u=f(u), &z=(x,y)\in\R^{m}\times \R^{N-m},\\
u(z)\in H^1(\R^N),
\end{cases}
\end{equation}
where $N\geq 2$, $1\leq m<N$.  The potential $V(z):=x_1^2+x_2^2+\cdots+x_m^2$ is called the partial
confinement if $m<N$ and the full confinement if $m=N$.

Since the 1970s, the nonlinear stationary Schr\"{o}dinger equation
\begin{equation}\label{eqs1.1-1}
-\Delta u+V(z)u=|u|^{p-2}u,\ \ z\in\R^{N},
\end{equation}
has been extensively studied, which was sparked by the seminal works of Berestycki and Lions \cite{BL1,BL2} and Strauss \cite{Strauss}.
Equation \eqref{eqs1.1-1} is of a variational nature, meaning that its solutions can be searched as critical points of the energy functional
$I_0: H^1(\R^N)\rightarrow \R$
  defined by
\[
I_0(u)=\frac{1}{2}\int_{\R^N}\big(|\nabla u|^2+V(z)u^2\big)dz-\frac{1}{p}\int_{\R^N} |u|^p dz.
\]
However,  the conventional variational techniques cannot be utilized in a conventional manner because of the absence of compactness, which is a result of the noncompactness of the group of translations acting on $\mathbb{R}^N$.
This challenge is characterized by the noncompactness of the embedding $H^1(\R^N)\hookrightarrow L^p(\R^N)$, which requires us to impose appropriate  conditions on $V(z)$.

If $V(z)$ is radially symmetric, then the restriction of $I_0$ to $H^1_r(\R^N)$, which consists of radially symmetric functions, restores compactness. As a result, the existence of a positive ground state solution and infinitely many (sign-changing) solutions to \eqref{eqs1.1-1} were established in \cite{BL1,BL2,BW2,CZ}.

Assuming that the potential is coercive, meaning that
$V\in C(\R^N,\R_+)$ and $V(z)\rightarrow +\infty$ as $|z|\rightarrow \infty$,
 it can be shown that the corresponding variational space is compactly embedded in
 $L^p(\R^N)$  for $2\leq p<2^*$.
  This allows us to prove the existence of a positive and negative solution, or infinitely many (sign-changing) solutions to \eqref{eqs1.1-1}, as shown in \cite{BW3,BW4,BPW,BLW,R}.

Nevertheless, if $V(z)$ is not radially symmetric or coercive, the compactness issue remains a challenge. Numerous studies have focused on this case, and various approaches have been developed to address the case
\begin{equation}\label{condition1}
0<V_\infty:=\lim\limits_{|z|\rightarrow\infty}V(z)<+\infty.
\end{equation}
We would like to highlight that the topological situation can be quite different depending on whether
 $V(z)$ approaches $V_\infty$ from below or from above.

In \cite{Lions1,Lions2}, Lions proved that when  $V(z)\leq V_\infty$, equation \eqref{eqs1.1-1} has a ground state solution at an energy level below the first level where the Palais-Smale condition fails, see also  \cite{R}. On the other hand, when  $V(z)>V_\infty$,
\eqref{eqs1.1-1} does not have a ground state solution. However, under a  suitable decay condition on $V(z)$, the existence of a positive solution can be established by using minimax methods and working at higher energy levels, see \cite{BL,BLions}.
Furthermore, Cerami et al. \cite{CPS, CMP} proved the existence of infinitely many positive and sign-changing solutions to equation \eqref{eqs1.1-1} under condition \eqref{condition1} and some suitable asymptotic behavior condition on $V(z)$.

After the pioneering paper \cite{CR}, where the potential $V(z)$ was assumed to be periodic, the existence of positive ground state solution and infinitely many geometrically distinct solutions were obtained in \cite{LWZ,SW} and the reference therein. We also refer to the papers \cite{AMPW,MPW,WY}, where the existence of infinitely many positive (or sign-changing) solutions for equation \eqref{eqs1.1-1} has been shown by adapting the Lyapunov-Schmidt method.

Recently, Bellazzini et al.\cite{BBJV}  investigated  the existence of normalized solutions
for  the nonlinear Schr\"odinger equation with a partial confinement
\begin{equation}\label{eqs1.5}
-\Delta u+(x_1^2+x_2^2)u-|u|^{p-2}u=\lambda u,\ \ z=(x_1,x_2,x_3)\in\R^{3},
\end{equation}
under the mass constraint
\begin{equation}\label{eqs1.6}
\int_{\R^3}|u|^2dz=c,
\end{equation}
where $u\in H^1(\mathbb{R}^3,\mathbb{C})$, $\frac{10}{3}<p<6$, $c>0$
is a constant and $(u,\lambda)$ is
a pair of unknowns with $\lambda$ being a Lagrange multiplier. They used a concentration compactness argument to prove the existence of a ground state solution for equations \eqref{eqs1.5}-\eqref{eqs1.6}, which is radially symmetric with respect to
$(x_1,x_2)$ and even decreasing with respect to $x_3$. Subsequently, Wei and Wu \cite{WW} showed that for sufficiently small
$c>0$, the problem \eqref{eqs1.5}-\eqref{eqs1.6} also admits a second positive solution, which is a mountain-pass solution.

This paper investigates multiple solutions  to  \eqref{eqs1.1}, including the existence of positive ground state solutions, saddle type nodal solutions, and the nonexistence of least energy sign-changing solutions together with symmetric properties. It is worth noting that such existence and nonexistence results for the nonlinear Schr\"odinger equation in the presence of a partial confinement have not been studied in the literature to our knowledge.

We assume that the nonlinear term $f\in {C}^1(\R,\R)$ satisfies the following hypotheses:

($f_1$)\,\,\,$f(s)=o(|s|)$ as $s\rightarrow 0$;

\vskip 0.06truein

($f_2$)\,\,\,$\lim\limits_{|s|\rightarrow +\infty}\frac{f(s)}{|s|^{2^*-1}}=0$, $2^*:=\frac{2N}{N-2}$
is the critical Sobolev exponent;

\vskip 0.06truein

($f_3$)\,\,\,there exists $\gamma>2$ such that $0<\gamma F(s)\leq f(s)s$ for $s\in\mathbb{R}$, where $F(s):=\int_0^s f(t)dt$;

\vskip 0.16truein

($f_4$)\,\,\,$\frac{f(s)}{|s|}$ is a strictly increasing function on $(-\infty,0)$ and $(0,+\infty)$.

\vskip 0.06truein


\vskip 0.06truein

To proceed with the presentation of our main results, we fix some definitions that will be useful in the sequel.
 We denote by $\mathcal{H}$ the
Sobolev space equipped with the inner product and norm
\[
(u,v):=\int_{\R^N}\big[\nabla u \cdot \nabla v +(x_1^2+x_2^2+\cdots+x_m^2)uv\big] dz,\ \ \|u\|:=(u,u)^{\frac{1}{2}}.
\]
Weak solutions of \eqref{eqs1.1} correspond to the critical points of the following energy functional $I : \mathcal{H}\rightarrow \R$ defined by
\begin{equation*}
I(u)=\frac{1}{2}\int_{\R^N}\big[|\nabla u|^2+(x_1^2+x_2^2+\cdots+x_m^2)u^2\big]dz-\int_{\R^N} F(u) dz.
\end{equation*}
Obviously, the functional $I$ belongs to $C^1(\mathcal{H},\R)$.
We say $u_0\in \mathcal{H}$
is a ground state solution of \eqref{eqs1.1}, if
\begin{equation}\label{eqs1.1-119}
I(u_0)=c:=\inf\limits_{\mathcal{N}}I(u),
\end{equation}
where
\begin{equation}\label{spaceN}
\mathcal{N}:=\big\{u\in \mathcal{H}\setminus\{0\} :\ \langle I'(u),u\rangle=0\big\}.
\end{equation}

\vskip 0.08truein

Our first main result is on the existence of ground state solution to the NLS with a partial confinement.

\begin{thm}\label{t1.1}
Assume $N\geq 2$, $1\leq m<N$ and ($f_1$)-($f_4$) hold. Then equation \eqref{eqs1.1} possesses a positive ground state solution $u\in \mathcal{H}\cap C_{loc}^2(\R^N)$. Moreover, there exists a point $y_0\in \R^{N-m}$ such that $u(x,y-y_0)$ is radially symmetric and decreasing with respect to $x$ and with respect to $y-y_0$, respectively.
\end{thm}

\begin{remark}
We would like to point out that the techniques of the Schwartz symmetrization and reflexion type arguments used in \cite{BBJV, Lopes} for establishing symmetry results for equation \eqref{eqs1.1} with a pure power nonlinearity are not applicable to the case of a
general  nonlinearity $f$. 
Instead, we adopt an alternative method and specifically depend on the method of moving planes that were outlined in \cite{CL-Book,ChenLiMa,H,GNN,LN}. This enables us to prove the symmetry property of
 $u(x, y-y_0)$.
\end{remark}

As equation \eqref{eqs1.1} only exhibits a partial confinement rather than the full confinement, we have developed a method of moving planes by establishing the Hopf lemma and the strong maximum principle for
$$
-\Delta u+\big(x_1^2+x_2^2+\cdots+x_m^2+c(z)\big)u(z)\geq 0
$$
in unbounded domains.

\vskip 0.05truein

The following result  addresses
 the nonexistence of least
energy sign-changing solution.
\begin{thm}\label{t1.1-1}
Assume $N\geq 2$, $1\leq m<N$ and ($f_1$)-($f_4$) hold. Then equation \eqref{eqs1.1} has no least energy sign-changing solutions.
\end{thm}


\begin{remark}


Theorem \ref{t1.1-1} is sharp and precise in indicating that there is no least energy sign-changing solution for equation \eqref{eqs1.1} with a partial confinement,  however, such a solution arises when the full confinement is implemented, see \cite{BLW}.

In the case of the  full confinement, $\mathcal{H}\hookrightarrow L^q(\R^N)$ is compact for $q\in (2,2^*)$,  which enables us to prove the existence of least energy sign-changing solution by solving the following
minimization problem
\begin{equation}\label{eqs-1.8}
m_0:=\inf\limits_{\mathcal{M}}I(u),
\end{equation}
where
\begin{equation}\label{eqs1.9}
\mathcal{M}:=\big\{u\in \mathcal{H},\ \ \text{and}\ \ u^+\not\equiv 0,\ \ u^-\not\equiv 0\ \ :\ \langle I'(u),u^+\rangle=\langle I'(u),u^-\rangle=0\big\}.
\end{equation}

Nevertheless, when dealing with the case of a partial confinement, equation \eqref{eqs1.1}
is translation invariant at least in one direction. Through an analysis based on energy comparison, we prove the nonexistence of least energy sign-changing solutions.
\end{remark}

\vskip 0.08truein

The third main result focuses on a distinct symmetric setting. We would like to clarify some notations. Let $G$ be a closed subgroup of the orthogonal group $O(m)\times I_{N-m}$ such that
\[
Fix(G)=\{0\},\ \ \text{where}\ \ Fix(G):=\big\{z\in \R^N:\ Az=z\ \text{for all}\ A\in G\big\},
\]
where  $O(m)$ is the group of orthogonal matrices in $\mathbb{R}^m$ and $I_{N-m}$  is the identity matrix in $\mathbb{R}^{N-m}$.
Furthermore, let $\tau:\ G\rightarrow\{-1,1\}$ be a group homomorphism. The pair $(G,\tau)$ gives rise to a group action of $G$ on $\mathcal{H}$, which is defined as follows:
\[
g\circ u(z):=\tau(g)u(g^{-1}z)\ \ \text{for each}\ g\in G, \ \text{and}\ z\in \R^N.
\]
The next result addresses solutions of equation \eqref{eqs1.1} in the context of an invariant subspace
\begin{equation}\label{eqHG}
\mathcal{H}_G :=\big\{u\in \mathcal{H}\ \big|\ g\circ u=u \ \ \text{for each}\ \ g\in G\big\}.
\end{equation}

\begin{theorem}\label{t1.20}
Assume $N\geq 3$, $1\leq m<N$ and ($f_1$)-($f_4$) hold. Let $G$, $\tau$ be as above, and suppose
that $\tau\equiv 1$ or that $G$ is finite. Then
\begin{equation}\label{eqs1.1-19}
c_G:=\inf\limits_{\mathcal{N}_G}I(u),
\end{equation}
is achieved, where
\begin{equation}\label{eqNG}
\mathcal{N}_G:=\big\{u\in \mathcal{H}_G\setminus\{0\} :\ \langle I'(u),u\rangle=0\big\}.
\end{equation}
Every minimizer $u\in \mathcal{N}_G$ of $I|_{\mathcal{N}_G}$ is a $G$-invariant solution of equation \eqref{eqs1.1} which obeys the variational characterization
\[
c_G=\inf\limits_{u\in \mathcal{H}_G}\sup\limits_{t>0}I(tu)=\inf\limits_{\gamma\in \Gamma}\max\limits_{t\in [0,1]}I(\gamma(t)),
\]
where the family of paths $\Gamma$ is defined by $
\Gamma:=\big\{\gamma\in C([0,1],\mathcal{H}_G)\ \big| \ \gamma(0)=0,\ \text{and}\ I(\gamma(1))<0\big\}$.
Moreover, the minimizer does not change sign if $\tau\equiv 1$, and up to translation in $y$ variable, every  minimizer $u$ is radially
symmetric with respect to $y$, that is $u(z)=u(x,|y|)$.
\end{theorem}

Saddle type nodal solutions are an important tool in many areas of science and engineering, and they have been used to solve a wide range of problems. They are especially useful for modeling systems that have sudden changes in behavior, and they can provide valuable insights into the dynamics of complex systems.
Inspired by \cite{LW,XW}, we explore the possibility of existence of saddle type nodal solutions.

To continue with the discussion on the existence of  saddle type nodal solutions, we define
\[
\mathcal{H}_{odd}:=\big\{u\in \mathcal{H}:\ \ \text{for almost every}\ (x_1,x',y)\in\R\times \R^{m-1}\times \R^{N-m}, \  \ u(-x_1,x',y)=-u(x_1,x',y)\big\}.
\]
For each positive integer $\ell$, denote
{\small
\begin{equation}\label{g-ell}
g_\ell:=
\begin{pmatrix}
     \cos\frac{2\pi}{\ell}  \quad -\sin\frac{2\pi}{\ell} & \quad\\
     \sin\frac{2\pi}{\ell}  \qquad \cos\frac{2\pi}{\ell} & \quad \\
     \quad  \quad  \quad & \mathbf{I}_{N-2}
\end{pmatrix}_{N\times N},
\end{equation}}
where $\mathbf{I}_{N-2}$ is the $(N-2)$-dimensional identity matrix. It is clear that $g_\ell\in O(N)$ is the counter-clockwise $\frac{2\pi}{\ell}$-rotation in $x_1x_2$-plane.

Let $\overline{G}$ be the cyclic group generated by $g_\ell$, that is
\[
\overline{G}:=\big\{g_\ell^{i}\ \big|\ i=0, 1,\cdots,\ell-1\big\}.
\]
The action of $G$ on a function $u: \R^{N}\to \R$ is defined by
\begin{equation*}
g\circ u(z)=u(g^{-1}z),\ \ \text{for each}\ g\in \overline{G}.
\end{equation*}

To investigate the equation  \eqref{eqs1.1}, we take a strategic approach by studying it in the $G$-invariant subspace
\[
\mathcal{H}_{\overline{G}} :=\big\{u\in \mathcal{H}_{odd}\ \big|\ g\circ u=u \ \ \text{for each}\ \ g\in \overline{G}\big\},
\]
which is a specific subset of the overall solution space that is preserved by the action of the symmetry group $G$.
We now define the minimization problem
\begin{equation*}
m_{\overline{G}}:=\inf\limits_{\mathcal{N}_{\overline{G}}}I(u),
\end{equation*}
with the constraint of Nehari manifold
\begin{equation*}
\mathcal{N}_{\overline{G}}:=\big\{u\in \mathcal{H}_{\overline{G}}\setminus\{0\} :\ \langle I'(u),u\rangle=0\big\}.
\end{equation*}

Subsequently, we  present our main results, which demonstrate the existence  of saddle type nodal solutions to \eqref{eqs1.1} and the symmetry of the minimizer.

%

\vskip 0.06truein

\begin{theorem}\label{t1.2}
Assume $N\geq 3$, $2\leq m<N$, and ($f_1$)-($f_4$) hold. Then for each integer $\ell\geq 1$, equation \eqref{eqs1.1} possesses one solution $u\in \mathcal{H}_{\overline{G}}\cap C_{loc}^2(\R^N)$ with exactly $2\ell$ nodal domains, and $I(u)=m_{\overline{G}}$. Furthermore, if
\begin{equation}\label{f-div}
f'(s)\leq C(s^{\sigma_{1}}+s^{\sigma_{2}}),\ \ \text{for each}\ s\geq 0,
\end{equation}
with constants $\sigma_{1},\sigma_{2}>0$,
then the minimizer $u$ is reflectionally symmetric with respect to the hyperplane containing $L_i$ for $i=0,1,\cdots,\ell-1$, where
\[
L_i=\Big\{(r\cos\theta_i,r\sin \theta_i,x'',y)\ \big|\ r\geq0,\ \theta_i=\big(\frac{i}{\ell}+\frac{1}{2\ell}+\frac{1}{2}\big)\pi
,\ x''\in\R^{m-2},\ \ y\in \R^{N-m}\Big\}.
\]
Up to a translation, this minimizer $u$ is radially
symmetric with respect to $x''=(x_3,x_4,\cdots,x_m)$ and with respect to $y$, respectively. Therefore, $u(z)=u(x_1,x_2,|x''|,|y|)$.
\end{theorem}

\vskip 0.05truein

The subsequent  result employs the higher-dimensional odd-symmetry to establish the existence of nodal solutions for equation \eqref{eqs1.1} with more intricate nodal domains. This result highlights the effectiveness of odd-symmetry in facilitating the discovery of solutions with intricate structures. The utilization of odd-symmetry enables the identification of nodal solutions that may not have been apparent through conventional methods. These solutions have significant implications in various fields, including science and engineering, where complex structures are often desired. Therefore, the discovery of such solutions can lead to advancements in technology and improve our understanding of complex systems.

To present  the result on the existence of nodal solutions, we begin by specifying certain notations.

Let $1\leq k\leq m$ and $\tau_i:\R^N\rightarrow \R^N$ be the reflection in the $i$-th coordinate direction, that
is, for each $z=(x_1,x_2,\cdots,x_m,y)\in \R^N$,
\[
\tau_i(z)=(x_1,x_2,\cdots,-x_i,\cdots,x_m,y)\ \ \text{for}\ i=1,2,\cdots,k.
\]

We study the equation \eqref{eqs1.1} in the
$k$-odd symmetric Sobolev space denoted by
\[
\mathcal{H}_k:=\Big\{u\in \mathcal{H}\ \big|\ u(\tau_i(z))=-u(z)\ \ \text{for each}\ 1\leq i\leq k\ \text{and}\ x\in \R^N\Big\},
\]
in which the functions are required to be antisymmetric with respect to the  variables  $x_1, x_2, \cdots, x_k.$

The minimization problem to be considered is
\begin{equation*}
m_k:=\inf\limits_{u\in\mathcal{N}_k}I(u),
\end{equation*}
 with the associated $k$-odd Nehari manifold defined as
\[
\mathcal{N}_k:=\big\{u\in \mathcal{H}_{k}\setminus\{0\} :\ \langle I'(u),u\rangle=0\big\}.
\]

%

\begin{theorem}\label{t1.3}
Assume $1\leq k\leq m<N$ and ($f_1$)-($f_4$) hold.
Then for each $1\leq k\leq m$,
equation \eqref{eqs1.1} has a weak solution $u\in \mathcal{H}_k\cap C_{loc}^2(\R^N)$ with exactly $2^k$ nodal domains such that $I(u)=m_k$.
Moreover, up to a translation, the minimizer $u$ is radially
symmetric with respect to $x^{(k)}:=(x_{k+1},x_{k+2},\cdots,x_m)$ and with respect to $y$, respectively. Therefore, $u(z)=u(x_1,x_2,\cdots,x_k,|x^{(k)}|,|y|)$.
\end{theorem}

\begin{remark}
When $k=1$ or $2$, Theorem \ref{t1.3} is equivalent to the case where $\ell=1$ or $2$ in Theorem \ref{t1.2}. Specifically, equation \eqref{eqs1.1} admits one odd solution with exactly two nodal domains or a saddle solution that changes its sign exactly four times, respectively.

It is clear that the nodal set of the solutions obtained from Theorem \ref{t1.2} and Theorem \ref{t1.3} contains at least one hyperplane.
These saddle type nodal solution demonstrates a different feature for the nonlinear Schr\"odinger equation with the partial confinement, since nontrivial solutions for the following nonlinear Schr\"odinger equation
\[
-\Delta u+u=f(u),\ \ z\in \R^N,
\]
can not vanish on a hyperplane \cite{EL}.
\end{remark}

The organization of the remaining parts of this article is as follows. In Section \ref{s2},  we investigate the existence and symmetry of positive ground state solutions to \eqref{eqs1.1} and establish Theorem \ref{t1.1}.
We also prove some preliminary results which  are necessary for our analysis, such as
the decay property, the existence and positiveness of the minimizer, the maximum principle at infinity, the Hopf lemma, and the strong maximum principle.
Furthermore, we investigate the nonexistence of the minimizer associated with least energy sign-changing solutions and then prove Theorem \ref{t1.1-1}. Finally, Section \ref{s4} is devoted to  the existence of saddle type nodal solutions and proving Theorem \ref{t1.20}--Theorem \ref{t1.3}.

\section{Existence of positive ground state solution}\label{s2}

In this section, we will present some lemmas that will be utilized in the proof of Theorem \ref{t1.1}. We begin by introducing a non-vanishing lemma of the minimizing sequence up to the translation with respect to $y\in\R^{N-m}$.

\begin{lemma}\label{lem2.1}(\cite[Lemma 3.4]{BBJV})
Assume that $\sup\limits_{n\geq 1}\|u_n\|<\infty$ and there exists $\varepsilon_0>0$ such that
\[
\lim\limits_{n\rightarrow\infty}\|u_n\|_{L^q(\R^N)}\geq \varepsilon_0
\]
for some $q\in (2,2^*)$. Then there exists a sequence $\{y_n\}\subset \R^{N-m}$, we have
\[
u_n(x,y-y_n)\rightharpoonup \bar{u}\neq 0\ \text{weakly in}\ \mathcal{H}.
\]
\end{lemma}

\vskip 0.12truein

\begin{lemma}\label{lem2.2}
Assume $N\geq 2$, $1\leq m<N$ and ($f_1$)-($f_4$) hold. Let $u$ be a weak solution of \eqref{eqs1.1}. Then
$u\in C_{loc}^{2,\gamma}(\R^N)$ for some $0<\gamma<1$ and satisfies
\[
\lim\limits_{|z|\rightarrow +\infty} u(z)=0.
\]
\end{lemma}
\noindent{\bf Proof.}
According to Theorem 1.11 in \cite{Li},
we know that $u$ belongs to $L^\infty(\R^N)\cap C_{loc}^{1,\alpha}(\R^N)$ for some $0<\alpha<1$, and
\[
\lim\limits_{|z|\rightarrow +\infty} u(z)=0.
\]
By using the classical Schauder regularity
estimates (see \cite[Theorem 4.6]{GT}), we then conclude $u\in C_{loc}^{2,\gamma}(\R^N)$ for some $0<\gamma<\alpha$.
\qed

\vskip 0.12truein

 To prove the existence of positive ground state solutions to \eqref{eqs1.1}, we need to show the positiveness  of the infimum of  the functional $I
 (u) $ on $
\mathcal{N}:=\big\{u\in \mathcal{H}\setminus\{0\} :\ \langle I'(u),u\rangle=0\big\}$.

\begin{lemma}\label{lem2.3}
Assume $N\geq 2$, $1\leq m<N$ and ($f_1$)-($f_4$) hold. Then $c>0$, where $c$ is defined in \eqref{eqs1.1-119}.
\end{lemma}
\noindent{\bf Proof.}
By ($f_1$)-($f_2$), for any small $\epsilon>0$, for some $2<p<2^{\ast}$, there exists a constant $C_{\epsilon}>0$ such that
\begin{equation}\label{eqs2.1-1}
|f(s)|\leq \epsilon(|s|+|s|^{2^{\ast}-1})+C_{\epsilon}|s|^{p-1},\ \ |F(s)|\leq \frac{\epsilon}{2}|s|^{2}+\frac{\epsilon}{2^{\ast}}|s|^{2^{\ast}}+\frac{C_{\epsilon}}{p}|s|^{p}.
\end{equation}
Then for each $u\in \mathcal{N}$, one has
\[
\|u\|^2=\int_{\R^N}\big[|\nabla u|^2+(x_1^2+\cdots+x_m^2)u^2\big]dz=\int_{\R^N}f(u)udz\leq
\epsilon\int_{\R^N}u^{2}dz+\epsilon\int_{\R^N}u^{2^{\ast}}dz+C_{\epsilon}\int_{\R^N}|u|^{p}dz.
\]
Thus,
\[
\|u\|^{2}\leq C(\|u\|^{2^{\ast}}+\|u\|^p),
\]
which implies that there exists a positive constant $\alpha>0$ such that
\[
\|u\|\geq \alpha>0,\ \ \text{for each}\ u\in\mathcal{N}.
\]
Then, for every $u\in \mathcal{N}$, we have
\begin{align*}
I(u)&=I(u)-\frac{1}{\gamma}\langle I'(u),u\rangle=
\big(\frac{1}{2}-\frac{1}{\gamma}\big)\|u\|^2+\frac{1}{\gamma}\int_{\mathbb{R}^{N}}\big(f(u)u-\gamma F(u)\big)dz\\
&\geq\big(\frac{1}{2}-\frac{1}{\gamma}\big)\|u\|^2\geq\big(\frac{1}{2}-\frac{1}{\gamma}\big)\alpha^2>0.
\end{align*}
Therefore
\begin{equation*}
c=\inf\limits_{\mathcal{N}}I(u)>0.
\end{equation*}\qed

By the general minimax principle  and Lemmas \ref{lem2.1}-\ref{lem2.3}, one can derive that  \eqref{eqs1.1} admits a ground state solution.

\begin{lemma}\label{lem2.4}
Assume $N\geq 2$, $1\leq m<N$ and ($f_1$)-($f_4$) hold. Then equation \eqref{eqs1.1} possesses a ground state solution.
\end{lemma}
\noindent{\bf Proof.}
Observing that the Nehari manifold $\mathcal{N}$ is of $C^1$ and  recalling that $c$ has a minimax
description in the form of the mountain pass level (see Theorem 4.2 of \cite{W1}):
\[
c=\inf\limits_{\gamma\in \Gamma}\max\limits_{t\in [0,1]}I(\gamma(t)),
\]
where the family of paths $\Gamma$ is defined by
\[
\Gamma:=\Big\{\gamma\in C([0,1],\mathcal{H})\ \big| \ \gamma(0)=0,\ \text{and}\ I(\gamma(1))<0\Big\}.
\]
Then it follows from the general minimax principle (Theorem 2.8 in \cite{W1}) that, there exists a sequence $\{u_n\}$ in $\mathcal{H}$ such that
\begin{equation}\label{eqs2.2}
I(u_n)\rightarrow c \ \ \text{and}\ \ I'(u_n)\rightarrow 0 \ \ \text{strongly in}\ \ \mathcal{H}^*,\,\, \mbox{as} \,\, n\to \infty,
\end{equation}
 where $\mathcal{H}^*$ denotes the dual space of $\mathcal{H}$. If $n$ is large enough, we can deduce from \eqref{eqs2.2} that
\begin{equation*}
\begin{aligned}
c+1+o_n(1)\|u_n\|&\geq
I(u_n)-\frac{1}{\gamma}\langle I'(u_n),u_n\rangle= \left(\frac{1}{2}-\frac{1}{\gamma}\right)\|u_n\|^2
+\frac{1}{\gamma}\int_{\mathbb{R}^{N}}\left(f(u_{n})u_{n}-\gamma F(u_{n})\right)dz\\
&\geq\left(\frac{1}{2}-\frac{1}{\gamma}\right)\|u_{n}\|^2,
\end{aligned}
\end{equation*}
which implies $\{u_n\}$ is bounded in $\mathcal{H}$. Thus,
by using \eqref{eqs2.2} again, we conclude
\begin{equation}\nonumber
0<c=\lim\limits_{n\rightarrow\infty}I(u_{n})=\lim\limits_{n\rightarrow\infty}\left(I(u_{n})-\frac{1}{2}\langle I'(u_{n}),u_{n}\rangle\right)=\lim\limits_{n\rightarrow\infty}\int_{\mathbb{R}^{N}}\left(\frac{1}{2}f(u_{n})u_{n}-F(u_{n})\right)dz,
\end{equation}
which  together with \eqref{eqs2.1-1} implies that
\[
\int_{\R^N}|u_n|^pdz\geq C>0, \ \ \text{for some}\ \ p\in (2,2^*).
\]
One can deduce from Lemma \ref{lem2.1}  that there exists a sequence $\{y_n\}\subset \R^{N-m}$ (if $m=N$, then we choose $y_n\equiv 0$) such that
\[
v_n:=u_n(x,y-y_n)\rightharpoonup u\neq 0\ \text{weakly in}\ \mathcal{H},\,\, \mbox{as}\,\, n\to \infty.
\]
Then $\{v_n\}$ is also a bounded $(PS)_c$ sequence of $I$, since equation \eqref{eqs1.1} is translation invariant with respect to the variable $y$. Thus, up to a subsequence, we can assume that
\begin{equation*}
\begin{cases}
v_n\rightharpoonup u,   &\text{weakly in $\mathcal{H}$}, \\
v_n\rightarrow u,  &\text{strongly in $L^q_{loc}(\R^N)$ for each $q\in (2,2^*)$}, \\
v_n\rightarrow u,  &\text{a.e. $\R^N$}
\end{cases}
\end{equation*}
as $n \to \infty.$ Obviously, $u$ is a nontrivial solution for equation \eqref{eqs1.1}. Therefore
\[
c\leq I(u)\leq \liminf\limits_{n\rightarrow \infty}\left[\left(\frac{1}{2}-\frac{1}{\gamma}\right)\|v_{n}\|^{2}+\int_{\mathbb{R}^{N}}\left(\frac{1}{\gamma}f(v_{n})v_{n}-F(v_{n})\right)dz\right]=\lim\limits_{n\rightarrow \infty}I(v_{n})=c,
\]
which implies $u\in \mathcal{N}$ and $I(u)=c$.
\qed

\vspace {0.35cm}

 As illustrated in the proof of Theorem \ref{t1.1}, the contradiction arguments are conducted during the process of moving planes. This implies that when working on an unbounded domain, it is essential to ensure that any negative minima obtained are always enclosed within a fixed ball for all solutions. To achieve this,  the following maximum principle at infinity is required.

\begin{lemma}\label{lem2.5}(Decay at infinity)
Assume that $\Omega$ is an unbounded domain in $\mathbb{R}^{N}$ and $w(z)$ is a classical solution of
\begin{equation}\label{equ2.3}
\begin{cases}
-\Delta w(z)+c(z)w(z)\geq 0, &z\in\Omega,\\
w(z)\geq 0,  &z\in\partial\Omega,\\
\lim\limits_{|z|\rightarrow\infty}w(z)=0, &z\in\Omega,
\end{cases}
\end{equation}
with
\begin{equation}\label{eq2}
\liminf\limits_{|z|\rightarrow\infty}c(z)\geq0.
\end{equation}
Then there exists a constant $R_{0}>0$ such that if
\[
\Omega\subset \R^N\setminus B_{R_0}(0),
\]
then
\begin{equation}\label{eq3}
w(z)\geq 0\ \ \text{for each}\ \ z\in \Omega.
\end{equation}
\end{lemma}

\begin{remark}
In the previous result of Decay at Infinity in
\cite[Theorem 6.2.3]{ChenLiMa}, the following conditions were required:
\begin{equation}\label{decay1}
\lim_{|z|\to \infty} w(z) |z|^s =0
\end{equation}
and
\begin{equation}\label{decay2}
c(z)>-\frac{s(n-2-s)}{|z|^2} \,\,\, \text{if}\ \, |z|\geq R_0,
\end{equation}
for some $0<s<n-2$.
However, Lemma \ref{lem2.5} allows us to generalize \eqref{decay1} to
\[
\lim_{|z|\to \infty}w(z)=0,
\]
provided that we assume a slightly stronger condition than \eqref{decay2} (see \eqref{eq2}), where we introduce some new auxiliary functions (see \eqref{eqnaf1} and \eqref{eqnaf2} below)
instead of $\frac{1}{|z|^s}$.
\end{remark}

\noindent{\bf Proof of Lemma \ref{lem2.5}.}
Since
\[
\liminf\limits_{|z|\rightarrow\infty}c(z)\geq 0,
\]
there exists a constant $R_{0}>0$ such that
\begin{equation}\label{equ2.8}
c(z)+\frac{1}{16}>0\ \ \text{for each}\ \ z\in \R^N\setminus B_{R_0}(0).
\end{equation}

Denote
$$
U_k=\left[2k\pi+\frac{\pi}{12}, 2(k+1)\pi-\frac{\pi}{12}\right] \cap \left\{x_1 \mid z \in \Omega \right\},\,\, k \in \mathbb{N}
$$
and
$$
V_k=\left(2(k+1)\pi-\frac{\pi}{12}, 2(k+1)\pi+\frac{\pi}{12}\right) \cap \left\{x_1 \mid z \in \Omega \right\},\,\, k \in \mathbb{N}.
$$
Obviously,
$$
\Omega= \mathop{\bigcup}\limits_{k \in \mathbb{N}} \left((U_k \times \mathbb{R}^{N-1}) \cap \Omega \right)\,\,  \cup \,\, \left((V_k \times \mathbb{R}^{N-1}) \cap \Omega \right).
$$
In order to prove \eqref{eq3}, it suffices to show
\begin{equation}\label{equ2.9}
w(z) \geq 0  \,\, \mbox{ in } \,\, (U_k \times \mathbb{R}^{N-1}) \cap \Omega
\end{equation}
and
\begin{equation}\label{equ2.10}
w(z) \geq 0  \,\, \mbox{ in } \,\, (V_k \times \mathbb{R}^{N-1}) \cap \Omega
\end{equation}
respectively for each $k \in \mathbb{N}.$

First we show \eqref{equ2.9} holds for each fixed $k \in \mathbb{N}.$ To do this, we introduce
\begin{equation}\label{eqnaf1}
\phi_k(z) =\phi_k(x_1)= \sin \left(\frac{x_1-2k\pi}{4}+\frac{\pi}{6}\right)\,\, \mbox{ for }\,\, z \in (U_k \times \mathbb{R}^{N-1}) \cap \Omega
\end{equation}
and define
$$
w_k(z)=\frac{w(z)}{\phi_k(z)}\,\, \mbox{ for }\,\, z \in (U_k \times \mathbb{R}^{N-1}) \cap \Omega
$$
for any $k\in\mathbb{N}$.

For some fixed $k \in \mathbb{N},$ suppose in the contrary of  \eqref{equ2.9}, then there exists a point  $p^k \in (U_k \times \mathbb{R}^{N-1}) \cap \Omega$  such that
$$
w(p^k)<0.
$$
Since $\phi_k(z)$ is positive, and bounded away from zero in $ (U_k \times \mathbb{R}^{N-1}) \cap \Omega$, then $\lim \limits_{|z|\to \infty} w_k(z)=0$ and
$$
\min_{\Omega} w_k(z) \leq w(p^k)<0.
$$
Therefore there exists a point $\bar p^k \in (U_k \times \mathbb{R}^{N-1}) \cap \Omega$ such that
$$
w_k(\bar p^k):=\min_{\Omega} w_k(z) \leq w(p^k)<0.
$$
At this point, we obtain
$$
w_k(\bar p^k)<0,\,\, -\Delta w_k(\bar p^k)\leq 0\,\, \mbox{ and }\,\, \nabla w_k(\bar p^k)=0,
$$
and by \eqref{equ2.3}, we have
\begin{equation*}
\begin{split}
0\geq -\Delta w_k(\bar p^k) &=\left[-\Delta w(\bar p^k)+\frac{\Delta \phi_k(\bar p^k)}{\phi_k(\bar p^k)}w(\bar p^k) \right]\frac{1}{\phi_k (\bar p^k)}\\
&\geq\left[-c(\bar p^k)w(\bar p^k)-\frac{1}{16}w(\bar p^k)\right] \frac{1}{\phi_k (\bar p^k)}\\
&=-\frac{w(\bar p^k)}{\phi_k (\bar p^k)} \left[c(\bar p^k)+\frac{1}{16} \right],
\end{split}
\end{equation*}
where we have used the fact that $\phi_k (\bar p^k) \in (0, 1)$ and $\frac{-\Delta \phi_k (\bar p^k)}{\phi_k (\bar p^k)}=\frac{1}{16}.$ It thus follows that
$$
c(\bar p^k)+\frac{1}{16} \leq 0.
$$
Now if $|\bar p^k|$ is sufficiently large, we derive a contradiction with
\eqref{equ2.8}, which implies \eqref{equ2.9}.

\vspace {.15cm}

Next we show \eqref{equ2.10} hold for each fixed $k \in \mathbb{N}.$
In this case, we introduce
\begin{equation}\label{eqnaf2}
\psi_k(z) =\psi_k(x_1)= \sin \left(x_1-2(k+1)\pi + \frac{\pi}{6}\right)\,\, \mbox{ for }\,\, z \in (V_k \times \mathbb{R}^{N-1}) \cap \Omega
\end{equation}
and define
$$
h_k(z)=\frac{w(z)}{\psi_k(z)}\,\, \mbox{ for }\,\, z \in (V_k \times \mathbb{R}^{N-1}) \cap \Omega,
$$
for any $k \in \mathbb{N}.$

By a similar argument as above, for some $k\in \mathbb{N}$, if there exists a point
$\bar q^k \in (V_k \times \mathbb{R}^{N-1}) \cap \Omega$ such that
$$
h_k(\bar q^k):=\min_{\Omega} h_k(z)<0.
$$
Then
\begin{equation*}
\begin{split}
0\geq -\Delta h_k(\bar q^k) &=\left[-\Delta w(\bar q^k)+\frac{\Delta \psi_k(\bar q^k)}{\psi_k(\bar q^k)}w(\bar q^k) \right]\frac{1}{\psi_k (\bar q^k)}\\
&\geq\left[-c(\bar q^k)w(\bar q^k)-w(\bar q^k)\right] \frac{1}{\psi_k (\bar q^k)}\\
&=-\frac{w(\bar q^k)}{\psi_k (\bar q^k)} \left[c(\bar q^k)+1 \right],
\end{split}
\end{equation*}
consequently
$$
c(\bar q^k)+ 1  \leq 0.
$$
Now if $|\bar q^k|$ is sufficiently large, we derive a contradiction with
\eqref{equ2.8}, which implies \eqref{equ2.10}.

Combining \eqref{equ2.9} with \eqref{equ2.10}, we arrive at \eqref{eq3}. This completes the proof of Lemma \ref{lem2.5}.
\qed

\vspace{0.15cm}

During the process of moving planes,
 it is necessary to establish a Hopf's lemma and a strong  maximum principle for the differential inequality \eqref{eq-1}. These principles play a pivotal role in facilitating the continuous movement of the planes towards the limiting position.

\begin{lemma}[Hopf lemma]\label{lem2.5-1}
	Let  $\Omega$ be a domain in $\mathbb{R}^N$ and  $u$ be a classical nonnegative solution of
\begin{equation}\label{eq-1}
-\Delta u + \big(x_1^2+x_2^2+\cdots+x_m^2+c(z)\big)u(z) \geq 0,\ \ z\in \Omega
\end{equation}
with bounded $c(z)$. Assume that there is a ball $B_{r_{0}}(z')$contained in $\Omega$ with a point $z^0 \in \partial \Omega \cap \partial B_{r_{0}}(z')$ and
\[
u(z)>u(z^0)=0, \forall z \in   B_{r_{0}}(z').
\]
Then
\[
\frac{\partial u}{\partial \nu} (z^0)<0.
\]
\end{lemma}

\begin{remark}
The condition $u(z)\geq 0$ is redundant if
$c(z)\equiv 0$, and it is sufficient to assume
$u(z^0)\leq 0$. Furthermore, \eqref{eq-1} can be extended to
\[
-\Delta u + \sum_{i=1}^N b_i(z) D_i u+\big(x_1^2+x_2^2+\cdots+x_m^2+c(z)\big)u(z) \geq 0,\ \ z\in \Omega,
\]
with bounded $b_i(z)$ for $ i=1, \cdots, N$.
\end{remark}

\noindent{\bf Proof of Lemma \ref{lem2.5-1}.}
Let $\lambda_1$ be the first positive eigenvalue of
\begin{equation*}
\begin{cases}
	-\Delta \phi (z) =\lambda_1 \phi(z), & z\in B_1(z^{0}),\\
	\phi(z)=0,& z\in \partial B_1(z^{0}),
	\end{cases}
	\end{equation*}
with the corresponding eigenfunction $\phi(z)>0$.

Define
$\psi(z)=\phi\left(\frac{z}{2r}\right),$ then
\begin{equation*}
\begin{cases}
-\Delta \psi (z) =\frac{\lambda_1}{ 4 r^2} \psi(z), &z\in B_{2r}(z^{0}),\\
\psi(z)=0, & z\in \partial B_{2r}(z^{0}).
	\end{cases}
	\end{equation*}
where $r$ is chosen to be satisfying
\begin{equation}\label{rrange}
0<r\leq \min\left\{ \frac{1}{2} \lambda_1^\frac{1}{2} \left(\mathop{\sup}\limits_{z\in \Omega} |c(z)|+1\right)^{-\frac{1}{2}},\frac{r_{0}}{2}\right\}.
\end{equation}

We introduce the auxiliary  functions $w$, $\tilde u$ and $v$ by defining
\[
w(z)=e^{-\alpha r^2}-e^{-\alpha |z-z'|^2},
\]
\[
\tilde u(z)=\frac{u(z)}{\psi(z)},
\]
and
\[
v(z)=\tilde u(z)+ \varepsilon w(z)
\]
in $D=B_{r}(z^0)\cap B_{r_{0}}(z'),$
where $\alpha$ and $\varepsilon$ are positive constants yet to be determined.

From \eqref{eq-1}, it is easy to deduce that
\begin{equation*}
\begin{aligned}
L \tilde u &:=-\Delta \tilde u-2 \nabla \tilde u \cdot \frac{\nabla \psi }{\psi}+ \left(\frac{-\Delta \psi}{\psi}
+x_1^2+x_2^2+\cdots+x_m^2+c(z)
\right)\tilde u(z)\\
&=-\Delta \tilde u-2 \nabla \tilde u \cdot \frac{\nabla \psi }{\psi}+ \left(\frac{\lambda_1}{4 r^2}
+x_1^2+x_2^2+\cdots+x_m^2+c(z)
\right)\tilde u(z)\\
&=-\Delta \tilde u+\sum_{i=1}^{N} b_i(z) D_i \tilde u +\tilde c(z) \tilde u \geq 0,
\end{aligned}
	\end{equation*}
where  $b_i(z)=-2\frac{D_i \psi(z)}{\psi(z)}$, $i=1, \cdots, N$ are bounded and $\tilde c(z)= \frac{\lambda_1}{4 r^2}
+x_1^2+x_2^2+\cdots+x_m^2+c(z)\geq 0$ due to \eqref{rrange}.

For any $w(z)$ in $D$, we  directly calculate
\begin{equation*}
\begin{aligned}
L w &=-\Delta w+\sum_{i=1}^{N} b_i(z) D_i w +\tilde c(z)w \\
&= e^{-\alpha |z-z'|^2}
\left[
4 \alpha^2 |z-z'|^2- 2\alpha \left(N-\sum_{i=1}^N b_i(z) (z_i-z'_{i})\right)
\right]
+\tilde c(z)\left(e^{-\alpha r^2}-e^{-\alpha |z-z'|^2}\right)\\
&\geq e^{-\alpha |z-z'|^2}
\left[
4 \alpha^2 |z-z'|^2- 2\alpha \left(N-\sum_{i=1}^N b_i(z) (z_i-z'_{i})\right)
\right],
\end{aligned}
	\end{equation*}
	here we have used the fact that $r\leq \frac{r_0}{2}\leq |z-z'|^2$ for any $z$ in $D$. Then
we can choose $\alpha $ sufficiently large such that $L w \geq 0$, it follows that
\begin{equation*}
Lv \geq 0 \,\, \mbox{ in } \,\, D,
\end{equation*}
then we can apply the weak maximum principle to conclude that
\begin{equation}\label{eq-2}
\min_{\partial D} v \leq \min_{\bar D} v.
\end{equation}
Therefore, we conclude that
 the minimum of $v$ on the boundary is actually attained at $z^0$:
\begin{equation}\label{eq-3}
v(z)\geq v(z^0),\,\, z \in \partial D.
\end{equation}

In fact, if $ z\in \partial D \cap \partial B_{r_0}(z')$, then
\[
v(z)=\tilde u(z) \geq \tilde u(z^0)=v(z^0)=0;
\]
if $z \in \partial D \cap B_{r_0}(z'),$ we can take $\varepsilon$ small enough such that
\[
v(z)=\tilde u(z)+\varepsilon w(z)\geq v(z^0)=0.
\]
This verifies \eqref{eq-3}.

Combining  \eqref{eq-2} with \eqref{eq-3} yields
\[
v(z) \geq v(z^0), \, \,  \text{for each}\ \ z\in D,
\]
it follows that
\[
\frac{\partial v}{\partial \nu}(z^0)\leq 0,
\]
noticing that
\[
\frac{\partial w}{\partial \nu}(z^0)> 0.
\]
Thus we arrive at the desired inequality
$$
\frac{\partial u}{\partial \nu}(z^0)<0.
$$
This completes the proof of Lemma \ref{lem2.5-1}.
\qed

\vspace{0.16cm}

\begin{lemma}[Strong maximum principle]\label{lem2.5-2}
	Let $\Omega $  be a domain in $\mathbb{R}^N$ with smooth boundary. Let $u$ be a classical nonnegative solution to
\begin{equation*}
\begin{cases}
-\Delta u + \big(x_1^2+x_2^2+\cdots+x_m^2+c(z)\big)u(z) \geq 0, & z\in \Omega,\\
u(z)=0, & z\in \partial \Omega,
\end{cases}
\end{equation*}
with bounded function $c(z)$. Then
	
(i) if $u$ vanishes at some point in $\Omega$, then $u \equiv 0$ in $\Omega$;
	
(ii) if $u \not\equiv 0$, then on $\partial \Omega$, the exterior normal derivative $\frac{\partial u}{\partial \nu}<0$.
\end{lemma}

\noindent{\bf Proof.}
(i) If $u$ vanishes at some point in $\Omega,$ but $u \not\equiv 0$ in $\Omega.$ Let
\[
\Omega_+=\big\{z\in \Omega \mid u(z)>0\big\}.
\]
Then $\Omega_+$ is an open set with $C^2$ boundary and
\[
u(z)=0,\ \ \text{for each} \, \, z\in \partial \Omega_+.
\]
Let $z^0$ be a point on $\partial \Omega_+$, but not on $\partial \Omega$. Then for sufficiently small $\rho >0$, there exists a ball $B_{\frac{\rho}{2}}(\bar z) \subset \Omega_+$ with $z^0$ as its boundary point. Now we use the Hopf lemma (Lemma \ref{lem2.5-1}) to conclude that the outward normal derivative at the boundary point  $z^0$ of $B_{\frac{\rho}{2}}(\bar z)$ satisfies 
\begin{equation}\label{eq-4}
\frac{\partial u}{\partial \nu}(z^0)<0.
\end{equation}

On the other hand, since $z^0$ is also a minimum of $u$ in the interior of $\Omega$, we must have
\[
\nabla u(z^0)=0.
\]
This contradicts \eqref{eq-4} and hence proves (i).

(ii) If $u \not\equiv 0$, we can derive the desired result by  (i) and the Hopf lemma (Lemma \ref{lem2.5-1}). This completes the proof of Lemma \ref{lem2.5-2}.
\qed

\vspace{0.36cm}


\vskip 0.1truein

\noindent{\bf Proof of Theorem \ref{t1.1}.}
Define
\[
I_{+}(u):=\frac{1}{2}\int_{\R^N}\big[|\nabla u|^2+(x_1^2+x_2^2+\cdots+x_m^2)u^2\big]dz-\int_{\R^N} F(u^+) dz,\ \ u\in \mathcal{H}.
\]
By a similar argument as Lemma \ref{lem2.4}, we conclude that equation \eqref{eqs1.1} admits a ground state solution $u\geq 0$. By Lemma \ref{lem2.2}, we know that $u\in C_{loc}^{2,\gamma}(\R^N)$ for some $0<\gamma<1$.
Moreover, it follows from Lemma \ref{lem2.5-2} that $u>0$. Therefore, equation \eqref{eqs1.1} admits a
positive ground state solution $u\in \mathcal{H}\cap C_{loc}^{2,\gamma}(\R^N)$.

We then proceed to investigate the symmetry properties of the positive ground state solution. To achieve this goal, we introduce the following notions.

For $\lambda\in \R$ and $z=(x,y)=(x_{1},x_{2},\cdots,x_{m},y)\in\mathbb{R}^{m}\times\mathbb{R}^{N-m}$, let
\[
{T}_{\lambda}:=\big\{\widetilde{z}=(\widetilde{x},\widetilde{y})\in \R^m\times\R^{N-m}\ |\ \widetilde{x}_{1}=\lambda\big\}
\]
be the moving planes,
\[
\Sigma_\lambda:=\big\{\widetilde{z}=(\widetilde{x},\widetilde{y})\in \R^m\times\R^{N-m}\  |\ \widetilde{x}_1<\lambda\big\},
\]
be the region to the left of the plane, and
$$
z^\lambda:=(2\lambda-x_{1},x_{2},\cdots,x_{m},y)
$$
be the reflection of $z$ about the plane $T_\lambda.$

Define
\[
{u}_\lambda(z):=u(z^\lambda)=u(2\lambda-x_{1},x_{2},\cdots,x_{m},y)
\]
and
\[
\Lambda:=\Big\{\lambda\in \R\ |\ u(z)<u(z^\lambda)\ \text{for all}\ z=(x,y)\in \R^m\times\R^{N-m}\ \text{with}\ x_1<\lambda \ \text{and}\ \frac{\partial u}{\partial x_1}>0\ \text{on}\ T_\lambda\Big\}.
\]

\vskip 0.12truein

Since $u$ goes to $0$ at $\infty$, there exists $R_1>R_0$ such that
\begin{equation}\label{eqs3.5}
\max\limits_{{B^c_{R_1}(0)}}u(z)<\min\limits_{\overline{B}_{R_0}(0)}u(z).
\end{equation}

\noindent{\bf Step 1.} $(-\infty,-R_1]\subset \Lambda$.

To compare the values of $u(z)$ with $u_\lambda(z),$ we denote
$$w_{\lambda}(z)=u(z^{\lambda})-u(z), \,\, z\in\Sigma_{\lambda}.$$ Then we have
\[
-\Delta u_{\lambda}+[(2\lambda-x_{1})^{2}+x_{2}^{2}+\cdots+x_{m}^{2}]u_{\lambda}=f(u_{\lambda}),\ z\in{\Sigma}_{\lambda}.
\]
Since $(2\lambda-x_{1})^{2}<x_{1}^{2}$ for $\lambda<0$, $z\in\Sigma_{\lambda}$, then
$$-\Delta u_{\lambda}+(x_{1}^{2}+x_{2}^{2}+\cdots+x_{m}^{2})u_{\lambda}\geq f(u_{\lambda}),\ z\in\Sigma_{\lambda}.$$
It follows that
\[
-\Delta w_{\lambda}+c(z)w_{\lambda}\geq0,
\]
where
\begin{equation}\label{eqcz}
c(z):=(x_1^2+x_2^2+\cdots+x_m^2)
-\int^{1}_{0}f'\big(\theta u_{\lambda}(z)+(1-\theta){u}(z)\big)d\theta.
\end{equation}
Then
\begin{equation*}
\begin{cases}
-\Delta w_\lambda+c(z)w_\lambda\geq 0\ &\text{in}\ \ \Sigma_\lambda,\\
w_\lambda=0   &\text{on}\ \ T_\lambda,\\
\lim\limits_{|z|\rightarrow\infty}w_\lambda(z)=0.
\end{cases}
\end{equation*}
Since $f'(0)=0$ and $\lim\limits_{|z|\rightarrow\infty}u(z)=0$,
we then deduce that
\[
\liminf\limits_{|z|\rightarrow\infty}c(z)\geq 0.
\]
Therefore, it follows from Lemma \ref{lem2.5} that
\[
w_\lambda(z)\geq 0\ \ \text{in}\ \Sigma_\lambda.
\]
On the other hand, if $\lambda\leq -R_1$, we deduce from \eqref{eqs3.5} that $w_{\lambda}>0$ on $\overline{B}_{R_0}(2\lambda,0,\cdots,0)\subset \Sigma_\lambda$. Hence, by using
Lemma \ref{lem2.5-2},
we conclude that
\begin{equation*}
\begin{cases}
w_\lambda>0\ \ &\text{in}\ \ \Sigma_\lambda,\\
\frac{\partial w_\lambda}{\partial x_1}<0\  \ &\text{on}\ \ T_\lambda.
\end{cases}
\end{equation*}
Thus, we prove that $(-\infty,-R_1]\subset \Lambda$.

\vskip 0.12truein

\noindent{\bf Step 2.} $\Lambda$ is open in $(-\infty,0)$.

Let $\lambda_0\in \Lambda\cap (-\infty,0)$. We will show that there exists an $\epsilon>0$ such that $(\lambda_0-\epsilon,\lambda_0+\epsilon)\subset \Lambda\cap (-\infty,0)$. We can assume $-R_{1}\leq\lambda_0<0$.

It follows from the assumption $\lambda_0\in \Lambda\cap (-\infty,0)$ that
\begin{equation}\label{eq:3.15}
w_{\lambda_{0}}(z)>0\ \ \text{in}\ \Sigma_{\lambda_0},\ \ \frac{\partial u}{\partial x_1}>0\  \ \text{on}\ T_{\lambda_0}.
\end{equation}
On $T_{\lambda_0}$, since $\frac{\partial u}{\partial x_1}>0$, we can find an $\epsilon_1\in (0,1)$ such that $\frac{\partial u}{\partial x_1}>0$ in
\[
\big\{z=(x_1,\cdots,x_{m},y)\in \R^m\times\R^{N-m}\ \big|\ \lambda_0-4\epsilon_1\leq x_1\leq \lambda_0+4\epsilon_1\ \ \text{and}\ \ |z'|\leq R_1+1\big\}
\]
with $z'=(x_2, \cdots, x_m, y)$.
Therefore we conclude that for any $\lambda\in (\lambda_0-\epsilon_1,\lambda_0+\epsilon_1)$, we have
\begin{equation}\label{eq:3.16}
\begin{cases}
w_{\lambda}(z)>0\ \ \text{in}\ \{z=(x,y)\in \overline{B}_{R_1+1}(0)\ |\ \lambda_0-2\epsilon_1\leq x_1<\lambda\},\\
\frac{\partial u}{\partial x_1}>0\  \ \text{on}\ T_{\lambda}\cap \overline{B}_{R_1+1}(0).
\end{cases}
\end{equation}
Let
\[
M:=2\max\left\{\left|\frac{\partial u}{\partial x_1}\right|\ \big|\ |x_{1}|\leq 2(R_1+1),\ |z'|\leq R_1+1\right\}
\]
and
\[
\delta:=\min\big\{u(z^{\lambda_0})-u(z)\ \big|\ |z'|\leq R_1+1,\ -(R_1+1)\leq x_1\leq \lambda_0-2\epsilon_1\big\}.
\]
Then \eqref{eq:3.15} implies that $\delta>0$, and hence
\begin{equation}\label{eq:3.17}
w_{\lambda}(z)>0\ \text{in}\ z\in \big\{z=(x,y)\in \overline{B}_{R_1+1}(0)\ |\ -(R_1+1)\leq x_1\leq \lambda_0-2\epsilon_1\big\}
\end{equation}
for any $\lambda\in (\lambda_0-\epsilon,\lambda_0+\epsilon)$ where
$\epsilon=\min\{\epsilon_1,\frac{\delta}{M},-\lambda_0\}$. Combining
\eqref{eq:3.16} with \eqref{eq:3.17}, for $\lambda\in (\lambda_0-\epsilon,\lambda_0+\epsilon)$, we obtain
\begin{equation}\label{eq:3.18}
\begin{cases}
w_{\lambda}(z)>0\ \ \text{in}\  \overline{B}_{R_1+1}(0)\cap\Sigma_\lambda,\\
\frac{\partial u}{\partial x_1}>0\  \ \text{on}\ T_{\lambda}\cap \overline{B}_{R_1+1}(0).
\end{cases}
\end{equation}
Now for any $\lambda\in (\lambda_0-\epsilon,\lambda_0+\epsilon)$, $w_{\lambda}\not\equiv 0$ in $\Sigma_\lambda\setminus \overline{B}_{R_1+1}(0)$ by \eqref{eq:3.18}. Then $w_{\lambda}$ satisfies
\begin{equation*}
\begin{cases}
-\Delta w_{\lambda}+c(z)w_{\lambda}\geq0\ \ \text{in}\ {\Sigma}_\lambda\setminus \overline{B}_{R_1+1}(0), \\
{w}_{\lambda}\geq 0\  \ \text{on}\ \partial\big({\Sigma}_\lambda\setminus\overline{B}_{R_1+1}(0)\big), \\
\lim\limits_{|z|\rightarrow\infty}w_{\lambda}(z)=0.
\end{cases}
\end{equation*}
Similarly as the proof in Step 1, for $\lambda\in (\lambda_0-\epsilon,\lambda_0+\epsilon)$, we have
\begin{equation*}
\begin{cases}
w_{\lambda}(z)> 0\ \ \ \text{in}\ \ {\Sigma}_\lambda\setminus\overline{B}_{R_1+1}(0), \\
\frac{\partial w_\lambda}{\partial x_1}<0\  \ \ \text{on}\ \ T_\lambda\setminus\overline{B}_{R_1+1}(0).
\end{cases}
\end{equation*}
which together  with \eqref{eq:3.18} yields
$(\lambda_0-\epsilon,\lambda_0+\epsilon)\subset {\Lambda}\cap(-\infty,0)$.
\vskip 0.12truein

\noindent{\bf Step 3.} We will show that either $\Lambda\cap (-\infty,0)=(-\infty,0)$ or $u(z^{\lambda_1})\equiv u(z)$ for some $\lambda_1\leq 0$.

Now we have shown that $\Lambda$ is open and contains all large negative $\lambda$ in $(-\infty,0)$. Let $(-\infty,\lambda_1)=\Lambda\cap (-\infty,0)$, so that $-R_{1}<\lambda_{1}\leq 0$. From the continuity of $u$, we have
\[
w_{\lambda_{1}}(z)=u(z^{\lambda_1})-u(z)\geq 0\ \text{in}\ \Sigma_{\lambda_1}.
\]
So we obtain
\begin{equation*}
\begin{cases}
-\Delta w_{\lambda_1}+c(z)w_{\lambda_{1}}\geq0\ \ \text{in}\ \Sigma_{\lambda_1}, \\
w_{\lambda_{1}}\geq 0\  \ \text{in}\ \ \Sigma_{\lambda_1},\ w_{\lambda_{1}}=0\  \ \text{on}\ T_{\lambda_1},\\
\lim\limits_{|z|\rightarrow\infty}w_{\lambda_{1}}(z)=0,
\end{cases}
\end{equation*}
where $c(z)$ is defined in \eqref{eqcz}.
By
Lemma \ref{lem2.5-2}, we have that either
\begin{equation}\label{eq:3.22}
w_{\lambda_{1}}\equiv 0\ \ \text{in}\ \Sigma_{\lambda_1}, \ \text{i.e.},\ u(z^{\lambda_1})\equiv u(z)\ \ \text{for}\ x_1<\lambda_1,
\end{equation}
or
\begin{equation}\label{eq:3.23}
\begin{cases}
w_{\lambda_{1}}>0\ \ &\text{in}\ \ {\Sigma}_{\lambda_1},\ \ \text{i.e.},\ u(z^{\lambda_1})>u(z) \ \ \text{in}\ \ {\Sigma}_{\lambda_1},\\
\frac{\partial w_{\lambda_1}}{\partial x_1}<0\ \ &\text{on}\ \ {T}_{\lambda_1}.
\end{cases}
\end{equation}
If \eqref{eq:3.23} happens, and if
$\lambda_1<0$, it would imply that $\lambda_1\in \Lambda\cap (-\infty,0)$, which is a contradiction. Hence $\lambda_1=0$, and \eqref{eq:3.23} becomes
\[
u(x_1,\cdots,x_{m},y)<u(-x_1,\cdots,x_{m},y)\ \ \text{for each}\ \ x_1<0,
\]
this completes the proof of Step 3.
\vskip 0.12truein

\noindent{\bf Step 4.} We show that $u$ must be symmetric in the $x_{1}$ direction about the hyperplane ${T}_{0}$. If \eqref{eq:3.22} occurs, then we have already shown that $u$
is symmetric in the $x_1$ direction about the hyperplane ${T}_{\lambda_1}$ and $\frac{\partial u}{\partial x_1}>0$ for each $x_1<\lambda_1$. On the other hand, if \eqref{eq:3.23} occurs with $\lambda_1=0$, we can repeat the previous Steps 2-3 for the positive $x_1$-direction for $u$ to conclude that either
\[
u^{\lambda_2}(z)\equiv u(z)\ \ \text{and}\ \ \frac{\partial u}{\partial x_1}<0\ \ \text{for each}\ \ x_1>\lambda_2\ \ \text{and}\ \
\lambda_2\geq 0,
\]
or
\[
u(x_1,\cdots,x_{m},y)>u(-x_1,\cdots,x_{m},y)\ \ \text{for each}\ \ x_1<0.
\]
However, the latter case and \eqref{eq:3.23} can not occur at the same time. Therefore $u$ must be symmetric in $x_1$ direction about some hyperplane $T_\lambda$ and strictly decreasing away from $T_\lambda$. Next since $u(z)$ is a solution of \eqref{eqs1.1}, then $u(z^{\lambda})=u(2\lambda-x_{1},x_{2},\cdots,x_{m},y)$ is also a solution, that is
\begin{equation}\nonumber
\left\{
\begin{aligned}
&-\Delta u(z^{\lambda})+(x_{1}^{2}+x_{2}^{2}+\cdots+x_{m}^{2})u(z^{\lambda})=f(u(z^{\lambda})),\\
&-\Delta u(z^{\lambda})+((2\lambda-x_{1})^{2}+x_{2}^{2}+\cdots+x_{m}^{2})u(z^{\lambda})=f(u(z^{\lambda})),
\end{aligned}
\right.
\end{equation}
which implies that $\lambda=0$.

Since equation \eqref{eqs1.1} is invariant under rotation with respect to $x$, we can choose any direction $(x_1,x_2,\cdots,x_m,0)\in \R^m\times\R^{N-m}$ as the $x_1$-direction. Therefore, we can conclude that $u(x,y)$ is symmetric about the hyperplane which contains $\{0\}\times\R^{N-m}$, perpendicular to that direction. Moreover, $u(x,y)$ is strictly decreasing away from that hyperplane. As a result, $u(x,y)$ must be radially symmetric and decreasing with respect to $x$.

Finally, repeating the previous Steps 1-4,
we conclude that, for some $y_0\in\R^{N-m}$, $u(x,y-y_0)$ must be radially symmetric and decreasing with respect to $y-y_0$.
\qed

\vskip 0.25truein


Next, we will demonstrate the nonexistence of least energy sign-changing solutions by proving that the minimizer of the functional $I$, restricted to $\mathcal{M}$, cannot be attained.

\vskip 0.15truein

\noindent{\bf Proof of Theorem \ref{t1.1-1}.}
Observe that if $u\in \mathcal{N}$, then $|u|\in \mathcal{N}$. 
Then
\begin{equation}\label{wei-2-2}
 c=\inf\big\{I(u) : u\in \mathcal{N} \ \text{and}\ u\geq 0\ a.e. \ \text{in}\ \R^N\big\},
\end{equation}
where $c$ and  $\mathcal{N}$ are  defined in \eqref{eqs1.1-119} and \eqref{spaceN}, respectively.

By the density argument, for each $u\in \mathcal{N}$ with $u\geq 0$, there exists a sequence $\{u_n\}\subset C_0^\infty(\R^N)$
with $u_n\geq 0$ such that $u_n\rightarrow u$ strongly in $\mathcal{H}$. Moreover, for each $n\geq 1$, there is a unique $t_n>0$
such that $t_n u_n\in \mathcal{N}$. It is easy to verify that $t_n\rightarrow 1$ as $n\rightarrow\infty$. This shows that
$t_n u_n\in \mathcal{N}\cap C_0^\infty(\R^N)$
with $t_n u_n\geq 0$ and $t_n u_n\rightarrow u$ strongly in $\mathcal{H}$ as $n\rightarrow\infty$.
Hence, it follows from \eqref{wei-2-2} that
\begin{equation*}
 c=\inf\big\{I(u) : u\in \mathcal{N}\cap C_0^\infty(\R^N) \ \text{and}\ u\geq 0\ a.e. \ \text{in}\ \R^N\big\}.
\end{equation*}
Therefore, for any $\varepsilon>0$, there exists a nonnegative function $u\in \mathcal{N}\cap \mathcal{C}_0^\infty(\R^N)$ such that
\[
I(u)\leq c+\varepsilon.
\]
Without loss of generality, we assume that the support of $u$ is contained in $B_R(0)$ for some $R>0$, namely
\[
supp\ u \subset B_R(0).
\]
For some $k\in \mathbb{N}$ and $k>R$, we define
\[
\widetilde{u}(z):=u(x,y_1+k,y')-u(x,y_1-k,y'), \ \ (x,y_1,y')\in \R^m\times\R\times\R^{N-m-1}.
\]
It is easy to see that $\widetilde{u}\in \mathcal{M}$ satisfies
\[
m_0\leq I(\widetilde{u})=I(u)+I(-u)\leq 2c+2\varepsilon,
\]
where
$m_0=\inf\limits_{\mathcal{M}}I(u)$ is defined in \eqref{eqs-1.8}.
Since $\varepsilon>0$ is arbitrary, we deduce that
\begin{equation}\label{wei-2-4}
m_0\leq 2c,
\end{equation}
where $\mathcal{M}$ is defined in \eqref{eqs1.9}. 

Now, suppose that there exists $w\in \mathcal{M}$ such that $I(w)=m_0$. Then
\[
m_0=I(w)=I(w^+)+I(w^-).
\]
By Theorem 4.3 in \cite{W1} and the maximum principle, one has $I(w^+)>c$, $I(w^-)>c$, i.e., $m_0>2c$,
which contradicts \eqref{wei-2-4}. Thus, the proof of Theorem \ref{t1.1-1} is completed.
\qed

\vspace {0.25cm}



\section{The $G$-invariant setting}\label{s4}

This section is devoted to the proof of  Theorem \ref{t1.20}--Theorem \ref{t1.3}. As part of the proof is similar to Theorem \ref{t1.1}, we omit some details for brevity.
We first show  the positiveness of infimum of the functional $I(u)$ in
$\mathcal{N}_G:=\big\{u\in \mathcal{H}_G\setminus\{0\} :\ \langle I'(u),u\rangle=0\big\}$.

\begin{lemma}\label{lem3.1}
Assume $N\geq 3$, $1\leq m< N$ and ($f_1$)-($f_4$) hold, then
$c_G>0$, where $c_G$ is defined in \eqref{eqs1.1-19}.
\end{lemma}
\noindent{\bf Proof.}
Since $\mathcal{N}_G=\mathcal{H}_G\cap \mathcal{N}$, we have $c_G\geq c$, where $\mathcal{N}_G,\,\mathcal{H}_G$ and $\mathcal{N}$ are defined in \eqref{eqNG}, \eqref{eqHG} and \eqref{spaceN} respectively.
Then $c_G>0$ follows  from Lemma \ref{lem2.3}. This completes the proof of Lemma \ref{lem3.1}.
\qed

\vskip 0.12truein

Next  we prove the existence of the solution to \eqref{eqs1.1} in $\mathcal{H}_G$.

\begin{lemma}\label{lem3.2}
Assume $N\geq 3$, $1\leq m<N$ and ($f_1$)-($f_4$) hold. Then equation \eqref{eqs1.1} possesses a solution $u\in \mathcal{H}_G$.
\end{lemma}
\noindent{\bf Proof.}
Observing that the Nehari manifold $\mathcal{N}_G$ is of $C^1$, we recall that $c_G$ has a minimax
description in the form of the mountain pass level (see Theorem 4.2 of \cite{W1}):
\[
c_G=\inf\limits_{u\in \mathcal{H}_G}\sup\limits_{t>0}I(tu)=\inf\limits_{\gamma\in \Gamma}\max\limits_{t\in [0,1]}I(\gamma(t)),
\]
where the family of paths $\Gamma$ is defined by
\[
\Gamma:=\Big\{\gamma\in C([0,1],\mathcal{H}_G)\ \big| \ \gamma(0)=0,\ \text{and}\ I(\gamma(1))<0\Big\}.
\]
By Lemma \ref{lem3.2}, one has $c_{G}>0$.
It follows from the general minimax principle (\cite[Theorem 2.8]{W1}) that, there exists a sequence $\{u_n\}$ in $\mathcal{H}_G$ such that
\begin{equation}\label{eqs3.4}
I(u_n)\rightarrow c_G \ \ \text{and}\ \ I'(u_n)\rightarrow 0 \ \ \text{strongly in}\ \ \mathcal{H}_G^*,\,\, \mbox{as}\,\, n\to \infty,
\end{equation}
where $\mathcal{H}_G^*$ denotes the dual space of $\mathcal{H}_G$. If $n$ is large enough, we deduce from \eqref{eqs3.4} that
\begin{equation*}
\begin{aligned}
c_G+1+o_n(1)\|u_n\|&\geq
I(u_n)-\frac{1}{\gamma}\langle I'(u_n),u_n\rangle= \left(\frac{1}{2}-\frac{1}{\gamma}\right)\|u_n\|^2+\frac{1}{\gamma}\int_{\mathbb{R}^{N}}\left(f(u_{n})u_{n}-\gamma F(u_{n})\right)dz\\
&\geq\left(\frac{1}{2}-\frac{1}{\gamma}\right)\|u_{n}\|^2,
\end{aligned}
\end{equation*}
which implies $\{u_n\}$ is bounded in $\mathcal{H}_G$. Thus,
by using \eqref{eqs3.4} again, we conclude
\begin{equation}\nonumber
0<m_{G}=\lim\limits_{n\rightarrow\infty}I(u_{n})=\lim\limits_{n\rightarrow\infty}\left(I(u_{n})-\frac{1}{2}\langle I'(u_{n}),u_{n}\rangle\right)=\lim\limits_{n\rightarrow\infty}\int_{\mathbb{R}^{N}}\left(\frac{1}{2}f(u_{n})u_{n}-F(u_{n})\right)dz,
\end{equation}
which togenther with \eqref{eqs2.1-1} yields
\[
\int_{\R^N}|u_n|^pdz\geq C>0, \ \ \text{for some}\ \ p\in (2,2^*).
\]
By Lemma \ref{lem2.1}, there exists a sequence $\{y_n\}\subset \R^{N-m}$ such that
\[
v_n:=u_n(x,y-y_n)\rightharpoonup \bar{u}\neq 0\ \text{weakly in}\ \mathcal{H}_G, \,\, \mbox{as}\,\, n\to \infty.
\]
 Then $\{v_n\}$ is also a bounded $(PS)_{c_G}$ sequence of $I$, since equation \eqref{eqs1.1} is translation invariant with respect to the variable $y$. Thus, up to a subsequence, we can assume that
\begin{equation*}
\begin{cases}
v_n\rightharpoonup \bar{u},   &\text{weakly in $\mathcal{H}_G$}, \\
v_n\rightarrow \bar{u},  &\text{strongly in $L^q_{loc}(\R^N)$ for each $q\in (2,2^*)$}, \\
v_n\rightarrow \bar{u},  &\text{a.e. $\R^N$},
\end{cases}
\end{equation*}
as $n \to \infty.$
Obviously, $\bar{u}$ is a nontrivial solution for equation \eqref{eqs1.1}. Hence
\[
c_{G}\leq I(\bar{u})\leq \liminf\limits_{n\rightarrow \infty}\left[\left(\frac{1}{2}-\frac{1}{\gamma}\right)\|v_{n}\|^{2}+\int_{\mathbb{R}^{N}}\left(\frac{1}{\gamma}f(v_{n})v_{n}-F(v_{n})\right)dz\right]=\liminf\limits_{n\rightarrow \infty}I(v_{n})=c_G,
\]
which implies $\bar{u}\in \mathcal{N}_G$ and $I(\bar{u})=c_G$.

By applying the symmetric criticality principle (\cite[Theorem 1.28]{W1} or \cite{Palais}), we conclude that $\bar{u}$ is a critical point of $I$ in $\mathcal{H}$. This completes the proof of Lemma \ref{lem3.2}.
\qed

\vspace {0.18cm}

\noindent{\bf Proof of Theorem \ref{t1.20}.}
It follows from Lemma \ref{lem3.2} that there exists a function $u\in \mathcal{H}_G$, such that $u\in \mathcal{N}_G$ is a minimizer of $I|_{\mathcal{N}_G}$. In the case of $\tau\equiv 1$, it is easy to see that $|u|\in \mathcal{N}_G$ is a minimizer of $I|_{\mathcal{N}_G}$ as well,
so it is a critical point of $I$ by the considerations above.
Moreover, it obeys the variational characterization stated in Theorem \ref{t1.20}.
By Lemma \ref{lem2.2}, $u\in C^{2,\gamma}_{loc}(\mathbb{R}^N)$, and 
$$-\Delta |u|+q(z)|u|=0, \,\, z\in \mathbb{R}^N$$
with $$q(z):=x_1^2+x_2^2+\cdots+x_m^2-\frac{f(u)}{u}\in L_{loc}^\infty(\R^N).$$ The strong maximum principle and the fact that $u\not\equiv 0$ therefore imply
$|u|>0$ on $\R^N$, which shows that $u$ does not change sign.

Therefore, we may assume $u>0$ in case of $\tau\equiv 1$. Then, following a similar argument as the proof of Theorem \ref{t1.1}, and  using the method of moving planes,
we conclude that, up to translation in $y$ variable, every  minimizer $u$ of $I|_{\mathcal{N}_G}$ is radially symmetric with respect to $y$, that is $u(z)=u(x,|y|)$.
\qed

\vskip 0.10truein

\noindent{\bf Proof of Theorem \ref{t1.2}.}
The existence of a saddle type solution $u\in \mathcal{H}_{\overline{G}}$ can be deduced from Theorem \ref{t1.20}, so we skip the details.

Now for any integer $\ell \geq 1$, we show that  $u$ has exactly $2\ell$ nodal domains.
Define the $\frac{\pi}{\ell}$-sectorial region by
\begin{equation*}
D:=\Big\{(r\cos\theta,r\sin\theta,x',y)\in\R^N\ \big|\ r\geq0,\ \theta \in\big[\frac{\pi}{2},\frac{\pi}{2}+\frac{\pi}{\ell}\big]
\ \ \text{and}\ \ x'\in \R^{m-2},\ \ y\in \R^{N-m} \Big\}.
\end{equation*}
We next show that $u$ has a constant sign in the interior of $D$ and $u$ is reflectionally symmetric with respect to the hyperplane stated in the Theorem \ref{t1.2}. We denote
\[
\mathcal{H}_0(D):=\Big\{u\in \mathcal{H}\ | \ u=0 \ \text{almost everywhere in}\ \R^N\setminus{int}(D)\Big\},
\]
and we define $I_\ell:\mathcal{H}_0(D)\rightarrow\R$ by
\[
I_\ell(v)=\ell\int_{D}\big[|\nabla v|^2+(x_1^2+x_2^2+\cdots+x_m^2)v^2\big]dz-2\ell\int_{D}F(v)dz.
\]
Obviously, for every $u\in \mathcal{H}_{\overline{G}}$, we have $\chi_D u\in \mathcal{H}_0(D)$ where $\chi_D$ is the characteristic function of $D$. Conversely, for each $v\in \mathcal{H}_0(D)$, we obtain
\[
\mathcal{S}(v)(z):=\mathcal{S}(v)(x,y)=\sum_{i=0}^{\ell-1}(\chi_D v)(g^i_\ell(x,y))-\sum_{i=0}^{\ell-1}(\chi_D v)(g^i_\ell\mathcal{R}(x,y))\in \mathcal{H}_{\overline{G}},
\]
where $g_\ell$ is given by \eqref{g-ell} and the reflection $\mathcal{R}$ is defined by
\begin{equation*}
\mathcal{R}=
\begin{pmatrix}
     -1\  \quad 0 & \quad\\
     \ \ 0  \quad\ 1 & \quad \\
     \quad  \quad  \quad & \mathbf{1}_{N-2}
\end{pmatrix}_{N\times N}.
\end{equation*}
Moreover, it can be checked directly that for any $u\in \mathcal{H}_{\overline{G}}$, $v\in \mathcal{H}_0(D)$, there holds
\[
u=\mathcal{S}(\chi_D u)\ \ \text{and}\ \ I_\ell(v)=I(\mathcal{S}(v)).
\]
These facts suggest that
\[
\inf\limits_{\mathcal{N}_D}I_\ell(u)
=\inf\limits_{\mathcal{N}_{\overline{G}}}I(u)=m_{\overline{G}},
\]
where
\[
\mathcal{N}_D=\big\{u\in \mathcal{H}_0(D)\setminus \{0\}\ |\ \langle I'_\ell(v),v\rangle=0\big\}.
\]

Let $u$ be a saddle type solution that achieves $m_G$, then $w=\chi_D u\in \mathcal{H}_0(D)$ and $I_\ell(w)=I(u)=m_{\overline{G}}$. Note that $|w|\in \mathcal{H}_0(D)$, we then derive by direct calculations that
\[
I_\ell(|w|)=I_\ell(w)=m_{\overline{G}}\ \ \text{and}\ \
\langle I_\ell'(|w|),|w|\rangle=\langle I_\ell'(w),w\rangle=0.
\]
We thus assume that $w\geq 0$ in $D$. On the other hand,
\[
\langle I'(\mathcal{S}(w)),\mathcal{S}(w)\rangle=\langle I_\ell'(w),w\rangle=0\ \ \text{and}\ \ I(\mathcal{S}(w))=I_\ell(w)=m_{\overline{G}}.
\]
This means that $\widetilde{u}:=\mathcal{S}(w)\in \mathcal{H}_{\overline{G}}$ also achieves $m_{\overline{G}}$, which implies $\widetilde{u}$ is a weak solution to the equation \eqref{eqs1.1-1}.

It follows from Lemma \ref{lem2.2} that $\widetilde{u}$ is of class $C^2_{loc}(\R^N)$ and
satisfies in the classical sense that
\[
-\Delta \widetilde{u}+(x_1^2+x_2^2+\cdots+x_m^2)\widetilde{u}\geq 0\ \ \text{in}\ D,
\]
which, together with Lemma \ref{lem2.5-2},
yields $\widetilde{u}>0$ everywhere in $int(D)$. The desired conclusion that $\widetilde{u}>0$ everywhere in $int(D)$ then follows from the definition of $\mathcal{S}$. In other words, $\widetilde{u}$ has exactly $2\ell$ nodal domains in the whole space $\R^N$.

Next, we investigate the symmetry property of the saddle solution $\widetilde{u}\in \mathcal{H}_{\overline{G}}$ by the method of rotating planes.
Recall that
\begin{equation*}
D=\Big\{(r \cos\theta, r\sin\theta,x'',y)\in\R^N\ |\ r\geq 0,\ \theta\in\left(\frac{\pi}{2},\frac{\pi}{2}+\frac{\pi}{\ell}\right),\ x''\in\R^{m-2},\ y\in \R^{N-m}\Big\}.
\end{equation*}
Fix $\theta\in (\frac{\pi}{2},\frac{\pi}{2}+\frac{\pi}{2\ell})$,
let $\Pi^\theta$ be the $N-1$ dimensional hyperplane defined by
\[
\Pi^\theta:=\Big\{(r \cos\theta, r\sin\theta,x'',y)\in\R^N\ |\ r\geq 0,\ x''\in\R^{m-2},\ y\in \R^{N-m}\Big\}.
\]
Define
\[
\Sigma_\theta:=\Big\{(r \cos\phi, r\sin\phi,x'',y)\in\R^N\ |\ r\geq 0,\ \phi\in (\frac{\pi}{2},\theta),\ x''\in\R^{m-2},\ y\in \R^{N-m}\Big\}.
\]
Let $z^\theta:=(x^\theta,y^\theta)$ denote the
reflection point of $z=(x,y)$ with respect to $\Pi^\theta$. Set
\[
\widetilde{u}_\theta(z):=\widetilde{u}(z^\theta)
\]
and
\[
w_\theta(z):=\widetilde{u}_\theta(z)-\widetilde{u}(z), \ \ \text{for each}\ z\in \Sigma_\theta.
\]
Clearly, $w_\theta(z)$ satisfies
\begin{equation}\label{eqs-system}
\begin{cases}
-\Delta w_\theta+(x_1^2+x_2^2+\cdots+x_m^2)w_\theta=c_\theta(z)w_\theta   &\text{in}\ \ \Sigma_\theta, \\
w_\theta=0   &\text{on}\ \ \Pi^\theta,\\
w_\theta>0   &\text{on}\ \ \partial\Sigma_\theta\setminus \Pi^\theta,
\end{cases}
\end{equation}
where
\[
c_\theta(z)=\frac{f(\widetilde{u}_\theta(z))-f(\widetilde{u}(z))}{\widetilde{u}_\theta(z)-\widetilde{u}(z)}=\int_0^1f'\big(t\widetilde{u}_\theta(z)+(1-t)\widetilde{u}(z)\big)dt.
\]
Consider the set
\[
\mathcal{A}^+:=\Big\{\theta\in (\frac{\pi}{2},\frac{\pi}{2}+\frac{\pi}{2\ell})\ |\ w_\theta(z)\geq 0\ \  \text{in}\ \ \Sigma_{\theta}\Big\},
\]
which is clearly a closed set in $(\frac{\pi}{2},\frac{\pi}{2}+\frac{\pi}{2\ell})$, since $w_\theta(z)$ is continuous with respect to $z$.

We claim that $\mathcal{A}^+$ is non-empty.
Observe first that $w_\theta^-:=\min\{w_\theta,0\}\in \mathcal{H}$. Moreover, for each $z\in\Sigma_\theta$
with $w_\theta<0$, it follows from \eqref{f-div} that
\begin{equation}\nonumber
\begin{aligned}
c_\theta(z)\leq C\int_0^1\big[\big(t\widetilde{u}_\theta(z)+(1-t)\widetilde{u}(z)\big)^{\sigma_{1}}+\big(t\widetilde{u}_\theta(z)+(1-t)\widetilde{u}(z)\big)^{\sigma_{2}}\big]dt\leq C\big[\widetilde{u}^{\sigma_{1}}(z)+\widetilde{u}^{\sigma_{2}}(z)\big].
\end{aligned}
\end{equation}
Also, the boundary conditions imply $w_\theta^-\equiv 0$ on
$\partial\Sigma_\theta$, and testing the equation
\eqref{eqs-system} by $w_\theta^-$ yields
\begin{align}\label{eqs-important}
\int_{\Sigma_\theta}|\nabla w_\theta^-|^2+(x_1^2+x_2^2+\cdots+x_m^2)|w_\theta^-|^2dz
&=\int_{\Sigma_\theta}c_\theta(z)(w_\theta^-)^2dz\nonumber\\
&\leq C\int_{\Sigma_\theta}\big[\widetilde{u}^{\sigma_{1}}(z)+\widetilde{u}^{\sigma_{2}}(z)\big](w_\theta^-)^2dz\nonumber\\
&\leq C_0\int_{\Sigma_\theta}(w_\theta^-)^2dz.
\end{align}
where $C_{0}=C\big(\|\widetilde{u}\|_{L^{\infty}(\Sigma_\theta)}^{\sigma_{1}}
+\|\widetilde{u}\|_{L^{\infty}(\Sigma_\theta)}^{\sigma_{2}}\big)$. Note that
$\widetilde{u}\in C^2_{loc}(\R^N)$, $\widetilde{u}=0$ on $\Pi^{\frac{\pi}{2}}$,
and $\lim\limits_{|z|\rightarrow\infty}\widetilde{u}(z)=0$, we deduce that
there exists $\delta>0$ such that
\[
C_{0}=C\big(\|\widetilde{u}\|_{L^{\infty}(\Sigma_\theta)}^{\sigma_{1}}
+\|\widetilde{u}\|_{L^{\infty}(\Sigma_\theta)}^{\sigma_{2}}\big)\leq \frac{1}{2}, \ \ \text{for each}\ \ \theta\in(\frac{\pi}{2},\frac{\pi}{2}+\delta).
\]
Consequently, $w_\theta^-\equiv 0$ provided that $\theta\in(\frac{\pi}{2},\frac{\pi}{2}+\delta)$
and this proves the claim.

Next, we claim that $\mathcal{A}^+$ is also open in $(\frac{\pi}{2},\frac{\pi}{2}+\frac{\pi}{2\ell})$.
To see this, let $\theta_0\in\mathcal{A}^+$. Since
$w_{\theta_0}\not\equiv 0$ by \eqref{eqs-system}, Lemma \ref{lem2.5-2} implies that $w_{\theta_0}>0$ in $\Sigma_{\theta_0}$. Note that, for each $p\in(2,2^*)$,
there exists $C_p>0$ such that
\[
\big(\int_{\R^N}|u|^p dz\big)^{\frac{2}{p}}\leq
C_p\int_{\R^N}\big[|\nabla u|^2+(x_1^2+x_2^2+\cdots+x_m^2)u^2\big] dz,\ \ \text{for each}\ u\in \mathcal{H}.
\]
Moreover, we may choose a compact set $D\subset \Sigma_{\theta_0}$ such that
\[
\|\widetilde{u}\|^{\sigma_{1}}_{L^{\tau_{1}}(\Sigma_{\theta_0}\setminus D)}+\|\widetilde{u}\|^{\sigma_{2}}_{L^{\tau_{2}}(\Sigma_{\theta_0}\setminus D)}\leq \frac{1}{2C_pC},
\]
where $C>0$ is the constant as in \eqref{eqs-important}.

On the other hand, by the continuity of the family $w_\theta$ with respect to $\theta$, there exists a
neighborhood $U\subset (\frac{\pi}{2},\frac{\pi}{2}+\frac{\pi}{2\ell})$ of $\theta_0$ with the property that
\[
w_\theta>0\ \ \text{and}\ \ \|\widetilde{u}\|^{\sigma_{1}}_{L^{\tau_{1}}(\Sigma_{\theta}\setminus D)}+\|\widetilde{u}\|^{\sigma_{2}}_{L^{\tau_{2}}(\Sigma_{\theta}\setminus D)}\leq \frac{1}{2C_pC}, \ \ \text{for each}\ \ \theta\in U,
\]
where $\tau_{i}=\frac{\sigma_{i}p}{p-2}$ for $i=1,\ 2$.
It follows from \eqref{eqs-important} and H\"older's inequality that
\begin{align*}
\big(\int_{\Sigma_\theta}|w_\theta^-|^pdz\big)^{\frac{2}{p}}
&\leq C_p\int_{\Sigma_\theta}\big[|\nabla w_\theta^-|^2+(x_1^2+x_2^2+\cdots+x_m^2)|w_\theta^-|^2\big]dz\nonumber\\
&=C_p\int_{\Sigma_\theta}c_\theta(z)(w_\theta^-)^2dz
\leq C_pC\int_{\Sigma_\theta}\big[u^{\sigma_{1}}(z)+u^{\sigma_{2}}(z)\big](w_\theta^-)^2dz\nonumber\\
&\leq C_pC\big(\|\widetilde{u}\|^{\sigma_{1}}_{L^{\tau_{1}}(\Sigma_{\theta}\setminus D)}+\|\widetilde{u}\|^{\sigma_{2}}_{L^{\tau_{2}}(\Sigma_{\theta}\setminus D)}\big)\big(\int_{\Sigma_\theta}|w_\theta^-|^pdz\big)^{\frac{2}{p}}\nonumber\\
&\leq
\frac{1}{2}\big(\int_{\Sigma_\theta}|w_\theta^-|^pdz\big)^{\frac{2}{p}},
 \ \ \text{for each}\ \ \theta\in U.
\end{align*}
Consequently, $w_\theta^-\equiv 0$ for $\theta\in U$ and this proves the claim.

\vskip 0.1truein

Since $\mathcal{A}^+$ is an open, closed and nonempty subset of $(\frac{\pi}{2},\frac{\pi}{2}+\frac{\pi}{2\ell})$, we conclude that
$\mathcal{A}^+=(\frac{\pi}{2},\frac{\pi}{2}+\frac{\pi}{2\ell})$.

Next, fix $\theta'\in(\frac{\pi}{2}+\frac{\pi}{2\ell},\frac{\pi}{2}+\frac{\pi}{\ell})$,
let $\Pi^{\theta'}$ be the $N-1$ dimensional hyperplane still defined by
\[
\Pi^{\theta'}:=\Big\{(r \cos\theta', r\sin\theta',x'',y)\in\R^N\ |\ r\geq 0,\ x''\in\R^{m-2},\ y\in \R^{N-m}\Big\}.
\]
Define
\[
\Sigma_{\theta'}:=\Big\{(r \cos\phi, r\sin\phi,x'',y)\in\R^N\ |\ r\geq 0,\ \phi\in (\theta',\frac{\pi}{2}+\frac{\pi}{\ell}),\ x''\in\R^{m-2},\ y\in \R^{N-m}\Big\}.
\]
Let $z^{\theta'}:=(x^{\theta'},y^{\theta'})$ denote the
reflection point of $z=(x,y)$ with respect to $\Pi^{\theta'}$. Set
\[
\widetilde{u}_{\theta'}(z):=\widetilde{u}(z^{\theta'})
\]
and
\[
w_{\theta'}(z):=\widetilde{u}_{\theta'}(z)-\widetilde{u}(z), \ \ \text{for each}\ z\in \Sigma_{\theta'}.
\]
Clearly, $w_{\theta'}(z)$ satisfies
\begin{equation*}
\begin{cases}
-\Delta w_{\theta'}+(x_1^2+x_2^2+\cdots+x_m^2)w_{\theta'}=c_{\theta'}(z)w_{\theta'}   &\text{in}\ \ \Sigma_\theta, \\
w_{\theta'}=0   &\text{on}\ \ \Pi^{\theta'},\\
w_{\theta'}>0   &\text{on}\ \ \partial\Sigma_{\theta'}\setminus \Pi^{\theta'},
\end{cases}
\end{equation*}
where
\[
c_{\theta'}(z)=\frac{f(\widetilde{u}_{\theta'}(z))-f(\widetilde{u}(z))}{\widetilde{u}_{\theta'}(z)-\widetilde{u}(z)}=\int_0^1f'\big(t\widetilde{u}_{\theta'}(z)+(1-t)\widetilde{u}(z)\big)dt.
\]
Consider the set
\[
\mathcal{A}^-:=\Big\{{\theta'}\in (\frac{\pi}{2}+\frac{\pi}{2\ell},\frac{\pi}{2}+\frac{\pi}{\ell})\ |\ w_{\theta'}(z)\geq 0\Big\},
\]
In the same manner, we can prove that
\[
\mathcal{A}^-:=\Big\{{\theta'}\in (\frac{\pi}{2}+\frac{\pi}{2\ell},\frac{\pi}{2}+\frac{\pi}{\ell})\ |\ w_{\theta'}({\color {red}z})\geq 0\Big\}=(\frac{\pi}{2}+\frac{\pi}{2\ell},\frac{\pi}{2}+\frac{\pi}{\ell}).
\]

Finally, a continuity argument also shows that $w_\theta\geq 0$ in $\Sigma_\theta$ for $\theta=\frac{\pi}{2}+\frac{\pi}{2\ell}$, and $w_{\theta'}\geq 0$ in $\Sigma_{\theta'}$ for $\theta'=\frac{\pi}{2}+\frac{\pi}{2\ell}$, which,
in particular, forces the symmetry of $\widetilde{u}$ with respect to the hyperplane
\[
L_0=\Big\{(r\cos\theta_0,r\sin \theta_0,x'',y)\ \big|\ r\geq0,\ \theta_0=\frac{\pi}{2}+\frac{\pi}{2\ell}
\ \ \text{and}\ \ x'\in \R^{m-2},\ y\in \R^{N-m}\Big\}.
\]

In view of our symmetric setting
again, $\widetilde{u}$ is reflectionally symmetric with respect to each hyperplane containing $L_i$ stated in
Theorem \ref{t1.2}.

Finally, by a similar argument as the proof Theorem \ref{t1.1},
up to translation in $y$ variable, we can prove that $\chi_D \widetilde{u}(x,y)$
is radially symmetric with respect to $x'':=(x_3,x_4,\cdots,x_m)$ and with respect to $y$, respectively. Therefore, we derive the desired symmetry property, that is, $\widetilde{u}(x,y)=u(x_1,x_2,|x''|,|y|)$. This 
completes the proof of Theorem \ref{t1.2}.
\qed

\vskip 0.15truein

\noindent{\bf Proof of Theorem \ref{t1.3}.}
We wish to highlight to the readers that the strategy employed for
the case that $\ell=1,2$ in Theorem \ref{t1.2} remains effective
for Theorem \ref{t1.3} as well.


Let $\tau_0$ be the identity map. By the definitions of the reflections, we easily see that
\[
\tau_i\tau_j=\tau_0\ \ \text{and}\ \ \tau_i\tau_j=\tau_j\tau_i,\ \
\text{for each}\ i,j=1,2,\cdots,k.
\]
Then for each $1\leq k\leq m$, the reflections $\tau_i$ for $i=1,2,\cdots,k$ generate a topological group $G_k$.
As a result, $G_k:=\langle\tau_i\ |\ i=1,2,\cdots,k\rangle$
is of order $2^k$. The action of $G_k$ on the space
$\mathcal{H}$ is defined for each $u\in\mathcal{H}$ by
\[
\tau_i\circ u(z)=-u(\tau_i (z))\ \ \text{for each}\ i=1,2,\cdots,k.
\]
Therefore, the $k$-odd symmetric Sobolev space $\mathcal{H}_k$ can be  rewritten as
\[
\mathcal{H}_k=\Big\{u\in\mathcal{H}\ |\ g\circ u=u \ \ \text{for each}\ g\in G_k\Big\}.
\]
Recall that
\[
m_k:=\inf\limits_{u\in \mathcal{N}_k}I(u),
\]
with the $k$-odd Nehari constraint $\mathcal{N}_k$ being defined by $\mathcal{N}_k=\mathcal{N}\cap\mathcal{H}_k$.
Therefore, the existence of a saddle type solution $u\in \mathcal{H}_k$ is a direct result of Theorem \ref{t1.20}.

Finally, we would like to point out that finding the minimizer of the
functional $I$ at the energy level $m_k$ on $\mathcal{N}_k$
is equivalent to searching for critical
point of the functional
\[
I_k(w)=2^{k-1}\int_{\widetilde{D}}\big[|\nabla w|^2+(x_1^2+x_2^2+\cdots+x_m^2)w^2\big]dz-2^{k}\int_{\widetilde{D}}F(w)dz.
\]
at the same energy level on the space $\mathcal{H}_0(\widetilde{D})$,
where $\widetilde{D}=\left([0,+\infty)\right)^k\times \R^{m-k}\times \R^{N-m}$.

By applying similar arguments as the proof of Theorem \ref{t1.2}, we may assume that $w=\chi_{\widetilde{D}} u\geq 0$ in $\widetilde{D}$. We then construct a weak solution $\bar{u}\in \mathcal{H}_k$ to the equation \eqref{eqs1.1} at the energy level $m_k$ by
\[
\bar{u}=\sum\limits_{\mathbf{k}}
(-1)^{\lambda}w(\tau_{\mathbf{k}}(z))\chi_{\widetilde{D}}(\tau_{\mathbf{k}}(z)),\ \ \text{with}\ \lambda=\sum_{j=1}^k k_j,
\]
where the multi-index $\mathbf{k}=(k_1,k_2,\cdots,k_k)$ with each
$k_j\in\{0,1\}$ and $j=1,2,\cdots,k$, $\tau_{\mathbf{k}}=\tau_1^{k_1}\tau_2^{k_2}\cdots \tau_k^{k_k}.$

It follows from Lemma \ref{lem2.2} that $\bar{u}$ is of class $C^2_{loc}(\R^N)$ and
satisfies in the classical sense that
\[
-\Delta \bar{u}+(x_1^2+x_2^2+\cdots+x_m^2)\bar{u}\geq 0\ \ \text{in}\ \widetilde{D},
\]
which, together with Lemma \ref{lem2.5-2}, yields $\bar{u}>0$ everywhere in $int(\widetilde{D})$.
In other words, $\bar{u}$ has exactly $2^k$ nodal domains in the whole space $\R^N$.

Finally, by a similar argument as the proof Theorem \ref{t1.1},
up to translation in $y$ variable, we can prove that $\chi_{\widetilde{D}} \bar{u}(x,y)$
is radially symmetric with respect to $x^{(k)}:=(x_{k+1},x_{k+2},\cdots,x_m)$ and with respect to $y$, respectively. Therefore, we obtain the desired symmetry property, that is, $u(x,y)=u(x_1,x_2,\cdots,x_k,|x^{(k)}|,|y|)$.
\qed
\vspace {0.45cm}

\noindent{\bf{Data availability statement.}}  Data sharing not applicable to this article as no datasets were generated or analysised during
the current study.
\smallskip\\

\noindent{\bf{Conflict of interest statement.}}  The authors declare that there have no conflict of interest.
\smallskip\\

\noindent{\bf{Acknowledgements.}} L. Shan and W. Shuai was partially supported by the national natural science foundation of China (No.12471107 and No. 12071170). L. Wu was partially supported by the national natural science foundation of China (No. 12401133 and No. 12031012).

\end{document}